\documentclass{article}
\usepackage{breakcites}
\usepackage{amsmath,amssymb,amsfonts,graphicx,bm,amsthm,mathrsfs}
\usepackage{authblk}
\usepackage[usenames]{xcolor}
\usepackage{soul}
\usepackage[colorlinks=true]{hyperref}
\definecolor{ultramarine}{rgb}{0 , 0, 0.3}
\definecolor{bulgarianrose}{rgb}{0.3, 0, 0}
\hypersetup{citecolor={ultramarine},
urlcolor={ultramarine},linkcolor={bulgarianrose}}
\usepackage{cleveref}
\usepackage{appendix}
\usepackage{pstricks}

\usepackage{enumitem}
\usepackage{latexsym}
\usepackage{stmaryrd}
\usepackage{bussproofs}
\usepackage{pstricks}
\usepackage{aliascnt}

\numberwithin{equation}{section} %adds section number before the equation label
\newtheorem{theorem}{Theorem}[section]
\newaliascnt{lemma}{theorem}  
\newtheorem{lemma}[lemma]{Lemma}  
\aliascntresetthe{lemma}  
  
\newaliascnt{corollary}{theorem}  
\newtheorem{corollary}[corollary]{Corollary}  
\aliascntresetthe{corollary}  
 
\newaliascnt{proposition}{theorem}  
  
\aliascntresetthe{proposition}  
 
\newtheoremstyle{exampstyle}
{\topsep} % Space above. The default value is \topsep
{\topsep} % Space below. The default value is \topsep
{} % Body font. For example you may add "\em" inside the braces.
{} % Indent amount
{\bfseries} % Theorem head font
{.} % Punctuation after theorem head
{.5em} % Space after theorem head
{} % Theorem head spec (can be left empty, meaning `normal')
\theoremstyle{exampstyle}% اینجا محیط تعریف شده در بالا خوانده میشود.
\newaliascnt{fact}{theorem}  
  
\aliascntresetthe{fact}  
 
\newaliascnt{claim}{theorem}  
  
\aliascntresetthe{claim}  
 
\newaliascnt{remark}{theorem}  
\newtheorem{remark}[remark]{Remark}  
\aliascntresetthe{remark}  
 
\newaliascnt{notation}{theorem}  
  
\aliascntresetthe{notation}  
 
\newaliascnt{example}{theorem}  
  
\aliascntresetthe{example}  
 
\newaliascnt{conjecture}{theorem}  
  
\aliascntresetthe{conjecture}  
 
\newtheorem{question}{Question}  
\newaliascnt{definition}{theorem}  
\newtheorem{definition}[definition]{Definition}  
\aliascntresetthe{definition}  
 
\newaliascnt{convention}{theorem}  
  
\aliascntresetthe{definition}

\newcommand{\Boxdot}{\,
                    \setlength{\unitlength}{1ex}
                    \begin{picture}(2,2)
                    \put(0,0){$\Box$}
                    \put(.55,.65){.}
                    \end{picture}
                    }  

\def\kcal{\mathcal{K}}

\def\nin{\not\in}

\def\HA{{\sf HA}}

\def\PA{{\sf PA}}

\def\K4{\hbox{\sf K4}{} }
\def\GL{\hbox{\sf GL}{} }

\def\R{\sqsubset}

\def\PL{{\sf PL}}
\def\lcal{\mathcal{L}}

\def\ar{\mid\!\sim}

\def\iph{{\sf{iPH}}}

\def\pce{\preccurlyeq}

\newcommand{\Ax}{\AxiomC}
\newcommand{\UI}{\UnaryInfC}
\newcommand{\BI}{\BinaryInfC}

\newcommand{\LLa}{\LeftLabel}
\newcommand{\RLa}{\RightLabel}
\newcommand{\DP}{\DisplayProof}

\newcommand{\lr }{\leftrightarrow}

\newcommand{\nat}{\mathbb{N}}

\def\nmodels{\not\models}
\newcommand{\myending}{\bibliography{\myref}%
\bibliographystyle{apalike}%
\end{document}%
}
\newcommand{\myheading}[1]{%
\begin{document}
\title{#1}
\author{
Mojtaba Mojtahedi\thanks{\url{http://mmojtahedi.ir}}\\
Ghent University\\
Department of Mathematics: Analysis, Logic and Discrete Mathematics}
\maketitle
}
\usepackage{mathtools}
\newcommand*{\upFrown}{\mathbin{\raisebox{0.9ex}{$\frown$}}}
\newcommand*{\upSmallFrown}{\mathbin{\raisebox{0.9ex}{$\smallfrown$}}}
\def\frown{\upSmallFrown}
\def\ILM{{\sf ILM}}
\def\SILM{\ILM^+}
\def\ILMS{{\ILM_\sigma}}
\def\SILMS{\ILM^+_\sigma}
\newcommand{\glpn}[1]{[{#1}]}
\newcommand{\glpnp}[1]{[{#1}]^+}
\newcommand{\glpp}[1]{\langle{#1}\rangle}
\def\R{\sqsubset}
\newcommand{\sub}[1]{{\sf sub}(#1)}
\def\sft{{\sf T}}
\def\tvdash{T\vdash}
\newcommand{\boolb}{{\sf Bool}(\Box)}
\def\natp{\nat^+}
\newcommand{\bool}{{\sf Bool}}
\newcommand{\mcal}{\mathcal{M}}
\newcommand{\ncal}{\mathcal{N}}
\def\emseq{\langle\, \rangle}
\newcommand{\invbox}[1]{\Box^{-1}#1}
\newcommand{\bas}[1]{{#1}^*}
\def\kfour{{\sf K4}}
\def\sfour{{\sf S4}}
\def\lcali{\lcal_\rhd}
\def\taut{{\sf Taut}}
\def\sfk{{\sf K}}
\def\four{{\sf 4}}
\def\lob{{\sf L\"ob}}
\def\cut{{\sf cut}}
\def\disj{{\sf disj}}
\def\conj{{\sf conj}}
\def\mont{{\sf Mont}}
\def\monts{{\sf Mont}_\sigma}
\def\ca{{\sf C_a}}
\def\boint{{\sf pint}}
\def\intcons{{\sf intc}}
\def\consint{{\sf cint}}
\def\nec{{\sf nec}}
\def\mods{{\sf mods}}
\def\boolb{{\sf Bool}(\Box)}
\def\sigb{{\Sigma_{\sf b}}}
\def\sigs{{\Sigma_\sigma}}
\def\sig{{\Sigma}}
\newcommand{\obox}[1]{{#1}^{\Box}}
\def\GL{{\sf GL}}
\def\GLS{\GL{\sf S}}
\def\Ca{{\sf C_a}}
\def\dhr{\rhd}
\def\rhdi{\rhd}
\def\Pre{{\sf Int}}
\def\Int{{\sf Int}}
\def\SPre{\Pre^+}
\def\lei{{\sf Lei}}
\newcommand{\myref}{bibliography.bib}
\newcommand{\mojref}{\bibliography{\myref}%
\bibliographystyle{apalike}%
}
\def\vdasht{\vdash_\sft}
\def\modelss{\models^{\!\sigma}}
\def\modelsa{\models}
\newcommand{\phis}[1]{\Phi^\sigma_{#1}}
\newcommand{\phia}[1]{\Phi_{#1}}
\newcommand{\succt}{\succ^t}
\newcommand{\prect}{\prec^t}
\def\pr{\mid\!\approx}
\newcommand{\pres}[2]{\mathrel{\setlength{\unitlength}{1ex}
                    \begin{picture}(2,2)
                    \put(0,0){$\pr$}
                    \put(-.1,1.3){\fontsize{4}{0} \selectfont ${#1}$}
                    \put(.2,-.8){\fontsize{4}{0} \selectfont ${#2}$}
                    \end{picture}}
                    }
\newcommand{\press}[2]{\mathrel{\setlength{\unitlength}{1ex}
                    \begin{picture}(2,2)
                    \put(0,0){$\pr$}
                    \put(-.1,1.3){\fontsize{4}{0} \selectfont ${#1}$}
                    \put(-.1,-.8){\fontsize{4}{0} \selectfont ${#2}$}
                    \end{picture}}
                    }                    
\newcommand{\modpres}[2]{\mathrel{\setlength{\unitlength}{1ex}
                    \begin{picture}(2,2)
                    \put(0,0){$\models$}
                    \put(-.1,1.3){\fontsize{4}{0} \selectfont ${#1}$}
                    \put(-.1,-.8){\fontsize{4}{0} \selectfont ${#2}$}
                    \end{picture}}
                    }                    
\newcommand{\nmodpres}[2]{\mathrel{\setlength{\unitlength}{1ex}
                    \begin{picture}(2,2)
                    \put(0,0){$\nmodels$}
                    \put(-.1,1.3){\fontsize{4}{0} \selectfont ${#1}$}
                    \put(-.1,-.8){\fontsize{4}{0} \selectfont ${#2}$}
                    \end{picture}}
                    }            
\def\V{\mathrel{V}}                              
\newcommand{\CPa}{{\sf CP_a}}
\def\lcalb{\lcal_\Box}
\def\ILMB{\ILM^\Box}
\newcommand{\ilmn}[1]{\ILM_{#1}}
\def\ILMSB{\ILMS^\Box}
\def\GLCa{\GL{\sf C_a}}
\newcommand{\glb}[3]{[{#3}]^{{#1}}_{#2}}
\newcommand{\adms}[2]{\mathrel{\ar}^{{\,{#1}}}_{\,{{#2}}}}
\newcommand{\rank}[1]{{\sf r}(#1)}
\newcommand{\trans}[2]{{#1}^{#2}}
\def\CL{{\sf CL}}
\def\sfK{{\sf K}}
\def\sce{\succcurlyeq}
\def\scomp{{\sf S}}
\def\psig{{\Sigma}}
\def\psigs{{\hat\psig}}
\def\NI{{\sf NI}}
\def\atom{{\sf atom}}
\def\pvar{{\sf var}}
\newcommand{\btra}[1]{{#1}^\Box}
\newcommand{\sbtra}[1]{{#1}^\boxdot}
\def\kfourp{{\kfour}{\sf p}}
\def\kfourpCa{{\kfourp}{\sf Ca}}
\def\GLp{{\GL}{\sf p}}
\def\GLpCa{\GLp{\sf Ca}}
\def\vdashs{\vdash_{\!\sigma}}
\def\sigz{\Sigma}
\newcommand{\hgt}[1]{h({#1})}
\def\Gama{\Delta}
\def\Gams{\Gamma_1}
\def\GLP{{\sf GLP}}
\newcommand{\boxn}[1]{[{#1}]}
\newcommand{\boxnp}[1]{[{#1}]^+}
\newcommand{\pre}[1]{{\sf pre}(#1)}
\def\fcal{\mathcal{F}}
\newcommand{\bsplit}[2]{\bgroup\def\arraystretch{1.2}\begin{tabular}{l}{#1}\\ {#2}\end{tabular} \egroup}
\def\sfour{{\sf S4}}
\def\bcal{\mathcal{B}}
\def\pcal{\mathcal{P}}
\newcommand{\lgc}[1]{{\sf T}_{#1}}
\newcommand{\lgcp}[1]{{\sf T}'_{#1}}
\newcommand{\plgc}[2]{\lgc{#1}^{#2}}
\newcommand{\blgc}[1]{{\bar{\sf T}}_{#1}}
\def\modelsp{\models^+}
\def\sfl{{\sf L}}
\def\pfrak{\mathfrak{P}}
\def\Rhat{\mathrel{\hat\R}}
\def\hatR{\mathrel{\hat\sqsupset}}
\def\Rone{\mathrel{\R^1}}
\def\Rsim{\mathrel{\sim_\R}}
\newcommand{\fel}[1]{\hat{#1}}
\newcommand{\iseg}[1]{\check{#1}}
\newcommand{\Bool}[1]{{\sf Bool}({#1})}
\def\vcal{\mathcal{V}}
\def\lcalo{\lcal_{\omega}}
\newcommand{\npred}[2]{{\sf pred}_{#1}(#2)}
\newcommand\pred{{\sf pred}}
\newcommand{\axioms}[1]{{\text{Axiom}}({#1})}
\newcommand{\rules}[1]{\text{Rule}({#1})}
\newcommand{\pnec}[1]{{\sf nec}_{#1}}
\newcommand{\plob}[1]{\text{L\"ob}_{#1}}
\def\lob{\plob{}}
\def\nec{\pnec{}}
\def\iglh{{\sf iGLH}}
\def\igl{{\sf iGL}}
\newcommand{\pmodelssymbol}{\setlength{\unitlength}{1ex}
                    \begin{picture}(2,2)
                    \put(-.3,0){$|$}
                    \put(0.1,.6){\rule{7pt}{0.3pt}} 
                    \put(0.1,0.2){\rule{7pt}{0.3pt}} 
                    \put(0.1,1){\rule{7pt}{0.3pt}}                                      
                    \end{picture}}                
\newcommand{\pmodels}{\mathrel{\pmodelssymbol}}
\newcommand{\npmodels}{\not\pmodels}   
\newcommand{\pmodelsp}{\pmodels^+}                 

\begin{document}
\title{Provability Models}
\author[1]{Mojtaba Mojtahedi\footnote{Email: mojtahedy@gmail.com.}}
\author[2]{Borja Sierra Miranda}
\affil[1]{Ghent University}
\affil[2]{University of Bern}
\maketitle

\begin{abstract}
In this paper, we study a new Kripke-style semantics for classical modal logic, 
named as provability models.
We study provability models for the propositional modal logics $\sfk$, $\kfour$, 
$\sfour$ $\GL$, $\GLP$ and the interpretability logic $\ILM$. 
Provability models combine features of Kripke models with the 
assignment of logics to individual worlds. 
Originally introduced in \cite{PLHA}, these models allowed the first author 
to establish arithmetical completeness for intuitionistic provability logic.

Interestingly, we show that the $\ILM$ is complete for the same provability models of 
$\GL$.
We improve provability models to predicative and decidable 
provability models in the case of $\GL$ 
and $\ILM$.
Furthermore, we prove a soundness and completeness of $\GLP$ for provability models.  
\end{abstract}
\tableofcontents

\section{Introduction}
Provability logic is a modal logic in which the modal 
operator is interpreted as \textit{provability} in a formal mathematical system. 
Several well-written texts address the topic of provability logic, most notably the two books 
\cite{Smorynski-Book,Boolos} and the two survey papers
\cite{ArtBekProv} and \cite{VisBek}.

In most studies of provability logic, the modal operator is typically interpreted as a provability predicate within some 
first-order arithmetical or set-theoretical system.
However, it is also possible to interpret $\Box A$ as 
$\sft\vdash A$ for a given propositional theory $\sft$.  
For instance, $\Box A\to \Box B$ could be interpreted as 
“$\sft\vdash A$ implies $\sft\vdash B$.” 
A precise definition of the $\sft$-provability 
interpretation for a propositional modal theory $\sft$
can be found in \Cref{ap-a}.

Although for any theory $\sft$ containing 
$\GL$\footnote{Here actually by this we mean that $\sft$ includes all theorems of $\GL$ and is closed under modus ponens and necessitation.},  one can show that $\GL$ is sound with respect to $\sft$-interpretations 
(\Cref{GL-sound-Lint}), no such $\sft$-interpretation yields completeness for $\GL$ 
(\Cref{GL-incompleteness}).
This failure, in contrast to $\PA$-interpretations,
stems primarily from the fact that the theory $\sft$ 
cannot simulate Kripke models. This contrasts with Peano Arithmetic, 
where Solovay was able to 
\cite{Solovay} exploit this 
capacity of $\PA$ to simulate Kripke models. 
Consequently, one cannot expect $\GL$ to serve as the provability logic of a single propositional theory. 

In this paper, we compensate for this limitation 
by explicitly incorporating the Kripke frame into the propositional 
theory in the following sense. 
We consider a standard Kripke model in which each world $w$ is augmented with a theory 
$\lgc w$. Naturally, we impose accessibility relationships between these theories based 
on the underlying accessibility relation (see \Cref{Mixed-sem}). 
Within this mixed model, the validity of boxed formulas $\Box A$ at a world $w$ 
is defined as the 
$\lgc u$-provability of $A$ in all worlds $u$ being accessible from $w$. 
This extension of propositional provability interpretations allows for both the 
soundness and completeness of $\GL$, 
as well as for many other modal logics.
We refer to such augmented models as 
\textit{provability models}.

Several concerns or objections regarding provability models are discussed in the following paragraphs. 
\\[2mm]
\textbf{Objection I:} 
\\
One might argue that provability models are merely a reformulation of standard Kripke models, based on the following reasoning: for a standard Kripke model, one could define $\lgc w$ as the set of all formulas valid at world $w$. Indeed, this is 
how completeness theorems are initially proved in this paper (see 
\Cref{lem1}). 
The concern is whether any provability model can be independent of standard Kripke models.  

The answer is affirmative, at least for modal logics of a provability-logic nature, namely 
$\GL$ and $\ILM$. In fact, for these logics, 
we show that for any finite irreflexive tree frame Kripke model, 
and any assignment of sets $X_w$ of formulas to the nodes $w$, one can construct a minimal 
$\lgc w\supseteq X_w$ such that the resulting structure constitutes a 
provability model (see \Cref{bases-gl}). In this sense, it is indeed possible to define 
new provability models independent from standard Kripke models. 

On the contrary, in case of the logic $\sfour$, the provability models correspond 
directly to the standard Kripke models (see \Cref{lem1,lem2}). 
\\[2mm]
\textbf{Objection II:}
\\
Is it possible to decide the validity of a formula at a given world in a provability model? 

In general, the answer is no. This is due to the possibility that the theories assigned to the nodes may be undecidable. 
Nevertheless, $\GL$ and $\ILM$ are complete for a class of decidable Provability models. Specifically, if one starts with a 
finite irreflexive tree frame Kripke model and assigns a finite set of formulas to each world, then it is possible to construct a decidable provability model 
(see \Cref{completeness-gl-finitary,completeness-ilm-finitary}).   
\\[2mm]
\textbf{Objection III:}
\\
What is the added value in studying provability models?

One important feature of provability models is their ability to absorb frame conditions as inference rules of the theories assigned to individual worlds. For instance, rather than requiring transitivity of the frame, one can require each theory to be closed under necessitation. 
Similarly, converse well-foundedness can be replaced by closure under L\"ob's rule. 

This feature was, in fact, the primary motivation for the first author to introduce 
provability models in the context of intuitionistic provability logic \cite{PLHA}.
Indeed, \cite{Iemhoff.Preservativity,Iemhoff,IemhoffT}
earlier introduced 
a form of Kripke semantics—known as Iemhoff semantics—for a candidate 
intuitionistic provability logic (the provability logic of Heyting Arithmetic $\HA$).
However, Iemhoff models are highly complex and hence ineffective for establishing
arithmetical completeness theorems. Their main drawback lies in their reliance 
on complicated frame properties. By contrast, the aforementioned feature of 
provability models allowed the first author of this paper to \textit{absorb} such frame conditions as 
inference rules of the theories assigned to the nodes, leading to finite and decidable 
provability models for intuitionistic provability logic. Thus, we have at least 
one concrete example, 
where provability models proved to be  particularly beneficial.

Another example is the provability models for the poly-modal provability logic $\GLP$.
$\GLP$ is a modal logic with infinitely many modal operators $\boxn n$, each intended to be 
a provability predicate for $\PA$ plus all true $\Pi_n$-sentences. 
Its arithmetical completeness was proved in \cite{Japaridze}. It is known that $\GLP$
is not complete for any Kripke frames. However, with the aid of provability models, and by 
absorbing some frame properties as logical properties, we are able to have soundness and 
completeness of $\GLP$ for provability models (see \cref{sec-GLP}).

\section{Definitions and elementary facts}
In this section we collect all basic definitions and simple facts about them. 
We separate them in different subsections, enabling the reader finding 
them more conveniently.
\subsection{Language}
In this paper, we deal with  three propositional modal languages. 
The first one, $\lcalb$, is with a unary modal operator $\Box$, 
together with boolean connectives $\to$ and $\bot$ and a set of atomic propositions 
$\atom$.  
The language $\lcali$, has a binary modal operator $\rhdi$ instead of $\Box$. 
The language $\lcalo $ for poly-modal provability logic 
$\GLP$, includes a unary modal operator 
$\boxn n $ for every $n\in\nat$.
We define the following notations:
\begin{itemize}
\item $\neg A:=A\to \bot$,
\item $A\vee B:=\neg A\to B$,
\item $A\wedge B:= \neg(\neg A\vee\neg B)$,
\item $A\lr B:= (A\to B)\wedge (B\to A)$,
\item  $\top:=\neg\bot$, 
\item $\lozenge A:= \neg\Box\neg A$,
%\item $A\rhdi B:=\neg B\dhr \neg A$,
\item $\Box A:=\neg A\rhdi \bot$ (in the language $\lcali$). Hence we may freely also use formulas of the language $\lcalb$, while we are working in the language $\lcali$.
\end{itemize}

We use $A,B,C,D,E,F,G$ for formulas in either of the mentioned languages. Also 
$p$ and $q$ are reserved for atomic propositions. 

A formula in the language $\lcalb$ ($\lcali$ or $\lcalo$) is called 
\textit{purely modal}, if it is boolean combination of formulas of the form 
$\Box B$ ($B_1\rhdi B_2$ or $\boxn n B$). In other words, a formula is purely modal, 
if every propositional variable is in the scope of some modal operator.
We will use the notation $\Boxdot A$ to mean the formula $A \wedge \Box A$.

If $\Gamma$ is a set of formulas we will write $\Box \Gamma$ to mean the set 
$\{\Box A \mid A \in \Gamma\}$.
Also, we will write $\Boxdot \Gamma$ to mean $\Gamma \cup \Box \Gamma$.

\subsection{Propositional theories and logics}\label{logicsec}
A  (propositional) \textit{theory}   $\sft$ (over the language $\lcal$) is a pair 
$(\axioms\sft,\rules\sft)$ with following properties:
\begin{itemize}
\item $\axioms\sft$ is a set of formulas  (in the language $\lcal$), called axioms 
of $\sft$.
\item $\rules\sft$ is a set of inference rules (for the language $\lcal$), called rules of $\sft$. 
The inference rules are considered as single-conclusion rules here, 
i.e.~they include a finite set of formulas as premise, and a single formula as 
conclusion.  
\end{itemize}
Then a \textit{logic} is just a theory which is closed under substitutions.
We have the following inference rules which occur in this paper often: (from left to right, they are modus ponens, necessitation and L\"ob's rule, respectively)
\begin{center}
\begin{tabular}{c c c}
\Ax{$A$}\Ax{$A\to B$}
\LLa{{\sf mp}}
\BI{$B$}
\DP
\hspace*{1cm}
&
\hspace*{1cm}
\Ax{$A$}
\LLa{{\sf nec}}
\UI{$\Box A$}
\DP
\hspace*{1cm}
&
\hspace*{1cm}
\Ax{$\Box A\to A$}
\LLa{L\"ob}
\UI{$A$}
\DP
\end{tabular}
\end{center}
Then we define the $\sft$-derivability relation 
$\vdash_\sft$ as follows: for a set $\Gamma\cup\{A\}$, we say $\Gamma\vdash_\sft A$
iff there is a finite sequence $(A_0,\ldots,A_n)$ such that every $A_i$ is either 
belonging to $\axioms\sft\cup \Gamma$ or it is derived by application of a rule 
in $\rules \sft$, with premises of the rule appeared in the sequence before $A_i$.  
When $\Gamma$ is empty, we may omit it in the notation $\Gamma\vdash_\sft A$ and simply 
write $\sft \vdash  A$.   
A \textit{theorem} of a theory $\sft$ is a formulas $A$ such that $\sft \vdash  A$.
%Then for the logic $\sft$ over the language $\lcal$, we defined 
%$$\bar\sft:=\{A\in\lcal:\ \  \vdash_\sft A\}.$$

For example, the Classical Logic, $\CL$, includes all Hilbert-style axioms of 
classical logic, together with  {\sf mp}.  
For the simplicity of notations, 
we use $\vdash$ for derivability in classical logic, i.e.~${\vdash}={\vdash_\CL}$. 

Also a theory $\sft$ is called \textit{classic}, if its rules of inferences includes 
modus ponens and all of the Hilbert-style axioms of classical logic are derivable 
in $\sft$.

We say that a theory $\sft$ includes $\sft'$, notation $\sft'\subset \sft$, if all theorems 
of $\sft'$ are also theorems of $\sft$. 

\subsection{Accessibility relation}\label{sec-order}

A relation $\R$ on a set $W$ is called \textit{converse well-founded}, if 
every non-empty subset $X\subseteq W$ has a maximal element, i.e.~some $w\in X$ such that there is no 
$u\in X$ with $w\R u$. Equivalently, it means that there is no infinite sequence 
$w_0\R w_1\R w_2\R\ldots$. 
All over this subsection we assume that $\R$ is a converse well-founded relation on $W$.

Given a relation $\R$, we use $\sqsupset$ for its inverse relation, i.e.~$w\sqsupset u$
iff $u\R w$.
Also define $\R^+$ as the transitive closure of $\R$. 
Similarly, $\sqsupset^+$ is the transitive closure of $\sqsupset$.
This means that 
$w\R^+ u$ if there is a sequence $w_0,\ldots,w_n$ with $n>0$ such that 
$w=w_0\R w_1\R\ldots\R w_n=u$.
Then we say that $(W,\R)$ is a \textit{tree}, if for every $w,u\R v$ we have either 
$w=u$ or $w\R^+ u$ or 
$u\R^+ w$.\footnote{Note that, with this definition, trees are not necessarily rooted.}

We say that a $w\in W$ is an \textit{immediate predecessor} of $u$, if $w\R u$ and 
there is no  $v\in W$ such that $w\R^+v\R^+ u$. In this situation we use notation 
$w\Rone u$. Note that by converse well-foundedness, 
any $\R$-accessible node, is also $\Rone$ accessible. Also if $(W,\R)$ is tree,
then the immediate predecessor should be unique. 
We denote  $\pred(w)$ as this unique immediate predecessor of $w$, 
in case of its existence.
We also define $w\Rsim u$ if either $w=u$ or there is some $v\Rone w,u$. 
Then we define $w\Rhat u$  iff   $w \Rsim v\R^+ u$ for some $v$.   Notice that $w\R^+ u$
implies $w\Rhat u$. 

\begin{lemma}\label{cwf-tree}
If  $(W,\R)$ is tree and converse well-founded, then  
$\Rhat$ is also converse well-founded.
\end{lemma}
%\begin{proof}
%Let $w_0\Rhat w_1\Rhat \ldots$ be an infinite sequence. Then we  build 
% another infinite sequence $w'_0\R^+ w'_1\R^+ \ldots$ violating the converse 
%well-foundedness of $(W,\R)$, a contradiction.
%
%Since $w_n\Rhat w_{n+1}$, either $w_n\R^+ w_{n+1}$ or 
%there are some $u,v\in W$ such that $u\Rone w_n,v$ and 
%$v\R w_{n+1}$. In the first case define  $w'_n:=w_n$  and 
%otherwise $w'_n:=u$. Then it is easy to observe that with the aid of quasi-finiteness, 
%we have $w'_0\R^+ w'_1\R^+\ldots $. 
%\end{proof}
\begin{proof}
    Note that since $\R$ is converse well-founded it is also irreflexive.
    In order to show that $\Rhat$ is converse well-founded assume we have a $\Rhat$-ascending chain $(w_i)_{i \in \mathbb{N}}$ and let us construct a $\R^+$-ascending sequence.
    We note that in a tree, $u \Rhat v$ means that either $u \R^+ v$ or that $\pred(u)$ exists and $\pred(u) \R^+v$.
    We define the sequence $(v_i)_{i \in \mathbb{N}}$ as
    \[
    v_i = \begin{cases}
        w_i &\text{if } w_i \R^+ w_{i+1}, \\
        \pred(w_i) &\text{if } \pred(w_i) \R^+ w_{i+1}.
    \end{cases}
    \]
    Note that this is well-defined since $w_i \Rhat w_{i+1}$ 
    and that, if it is not the case of $w_i \R^+ w_{i+1}$, then $\pred(w_i)$ exists.
    Also, by definition, we have that $v_i \R^+ w_{i+1}$.
    Since $v_i = w_i$ or $v_i = \pred(w_i)$ it follows from 
    $v_i \R^+ w_{i+1}$ that $v_i \sqsubseteq^+ v_{i+1}$.
    Consider the sequence $(v_{2i})_{i \in \mathbb{N}}$, 
    we will show that it is $\R^+$-ascending.
    We already know that $v_{2i} \sqsubseteq^+ v_{2i+1} \sqsubseteq^+  v_{2(i+1)}$, 
    so $v_{2i} \sqsubseteq^+ v_{2(i+1)}$, so to prove that $v_{2i} \R^+ v_{2(i+1)}$ 
    it suffices to show that $v_{2i} \neq v_{2(i+1)}$.
    Suppose that $v_{2i} = v_{2(i+1)}$, then the fact that $v_{i} = w_i$ or 
    $v_{i} = \pred(w_i)$ for any $i$, 
    together with $w_{2i} \R^+ w_{2i+1} \R^+ w_{2(i+1)}$,
    implies a contradiction. 
\end{proof}

\subsection{Finite sequences}\label{Fin-seq}
Since  in this paper, we need to manipulate 
finite sequences, it is 
more convenient to have some notations for 
the manipulation of finite sequences. We use $\sigma$,
$\tau$ and $\eta$ for finite sequences. 
We define:
\begin{itemize}
\item $\sigma*\tau$ as the concatenation of finite sequences $\sigma$ and $\tau$, 
\item $\fel \sigma$ as the final element of 
$\sigma$,
\item $\iseg \sigma$ as the sequence $\sigma$ without its final element.
\end{itemize}
\subsection{Kripke models}\label{def-KP}
A Kripke model, is a tuple $\kcal=(W,\R,\V)$ with following properties:
\begin{itemize}
\item $W$ is  a non-empty set of worlds. 
\item $\R$ is a binary relation on $W$.
\item ${\V}\subseteq W\times\atom$.
\end{itemize}
Then we define $\kcal,w\models B$ by induction on $B$ as follows. 
For atomic proposition $p$ define $\kcal,w\models p$ iff $w\V p$. 
If $B$ is having an outmost 
boolean connective, then define its valuation as in the classical valuation.\footnote{In particular, we never have $\mathcal{K}, w \models \bot$.}
In other  words, $\kcal,w\models$ commutes with boolean connectives. 
Finally $\kcal,w\models \Box B$ iff for every $u\sqsupset w$ we have $\kcal,u\models B$.
A frame condition for $\kcal$, always refers to $(W,\R)$.

Finally, we define
\[
\kcal,w \models^+ A \quad \text{iff} \quad \exists\, u \sqsubset w\ \forall\, v \sqsupset^+ u\ (\kcal, v \models A).
\]

\section{Provability models for the language $\lcalb$}\label{sec2}

In this section we first define provability models for the language $\lcalb$.
Then we prove the soundness and completeness of the modal logics 
$\sfk$, $\kfour$, $\sfour$ and $\GL$ for provability models. 
So let us begin with the definition of provability models:

\begin{definition}\label{Mixed-sem}
A \textit{provability pre-model} (pre-model in short) is a tuple 
$\pcal=(W,\R,\{\lgc w\}_{w\in W},\V)$ with following properties:
\begin{itemize}
\item $\kcal_\pcal:=(W,\R,\V)$ is a Kripke model (as defined in \cref{def-KP}).
%\item $\R$ is a relation on $W$.
%\item $W^\R$ is the set of all $\R$-accessible nodes, 
%i.e.,~$W^\R:=\{w\in W:\exists\, u\in W\ (u\R w)\}$.
\item  For any $\R$-accessible $w\in W$,
%$w\in W^\R$,
$\lgc w$ is a theory.   If no confusion is likely, we may use $\vdash_w$ for 
$\lgc w$-derivability relation $\vdash_{\lgc w}$.
%\item $\V$ is a valuation for variables, i.e.~a binary relation on $W$ and propositional variables.
\end{itemize}
\textit{Notation:} We use $W^\R$ for the set of all $\R$-accessible nodes in $W$.
Then we define $\pcal,w\pmodels A$ by induction on the complexity of $A$ 
such that it commutes with the boolean operators and 
$$\pcal,w\pmodels \Box A \quad\text{iff}\quad \forall u\sqsupset w \ \vdash_{u} A.$$
Moreover, we define 
$$\pcal,w\pmodelsp A \quad \text{iff}\quad \exists\, u\sqsubset w\ \forall\, 
v\sqsupset^+ u\ (\pcal,v\pmodels A).$$
Notice that if $w$ is not $\R$-accessible, then $\pcal,w\pmodelsp A$ does not hold for any
$A$.\\
$\pcal$ is called \textit{provability model}, if it is a pre-model with following property: 
\begin{itemize}
\item \textit{Modal completeness}: 
	For every purely modal $A$, if $\pcal,w\pmodelsp A$ then $\vdash_w A$.
%\item \textit{Classicality}: Each $\lgc w$ is classic, i.e.~they have modus ponens and all
%classical tautologies are derivable.
\end{itemize}
Define $\pcal\pmodels A$ if $\pcal,w\pmodels A$ for every $w\in W$. 
Moreover, for a class $\pfrak$ of pre-provability  models, 
$\pfrak\pmodels A$ if for every $\pcal\in\pfrak$ we have $\pcal\pmodels A$.  
Also $\pfrak\pmodels\Gamma$ for a set of formulas, 
if $\pfrak\pmodels A$ for every $A\in\Gamma$.
Moreover,   $\Gamma\pmodels_\pfrak A$ 
iff  $\pcal,w\pmodels \Gamma$ implies $\pcal,w\pmodels A$ for every $\pcal\in\pfrak$ and 
any world $w$ in $\pcal$. 
This notion  of semantical consequence relation, 
 is usually referred as \textit{local} in the 
literature (see \cite{Shillito}). This is in contrast with the \textit{global} version which is defined as ``$\pfrak\pmodels \Gamma$ implies $\pfrak\pmodels A$". Nevertheless, in this paper, we always consider the local version.\footnote{We thank Ian Shillito for bringing this subtle distinction to our attention.}

%\begin{itemize}
%\item $\pcal,w\pmodels A$,
%\item if $w\in W^\R$ then $\vdash_w A$.
%\end{itemize}

We say that a theory $\sft$ is \textit{sound} for a  provability model 
$\pcal$, if for every  $ \sft \vdash A$  we have $\pcal\pmodels A$.
Then we say that a theory is sound for a class $\pfrak$ of provability models, 
if it is sound 
for every member of it. We say that $\sft$ is \textit{complete} for $\pfrak$ if 
the following holds:  $\pfrak\pmodels A$ implies  $\sft\vdash A$.
Moreover, we say that $\sft$ is strongly complete for $\pfrak$, if 
for any set $\Gamma$ of formulas, $\Gamma\pmodels_\pfrak A$ implies $\Gamma\vdash_\sft A$.

Finally two provability (pre-) models are being said to be \textit{l-isomorphic}, if they share the same Kripke models (Same frame and evaluation on atomics) and the theories attached to the worlds are having the same set of derivable formulas. Note that, theories are having a set of axioms and a set of inference rules. It can be the case that you have different theories with the same sets of derivable theorems. 

A \textit{frame property} for a provability model, is a property that applies to 
$(W,\R)$. For example, the frame properties 
like reflexivity, transitivity or converse well-foundedness, means that $(W,\R)$ is so. 

On the other hand, we also have \textit{logical} properties of the provability model, 
which are properties that corresponds to all $\lgc w$'s for $w\in W^\R$. 
For example, a provability model with necessitation, means a provability model that 
all the theories $\lgc w$ are having necessitation as part of their inference rules.
Also we say that a provability model is \textit{classical}, if all of $\lgc w$'s are 
classical\footnote{Remember that this means that $\lgc w$ includes modus ponens as 
its rules of inferences, and includes all Hilbert-axioms of classical logic as its derivable formulas.}.
\\[3mm]
\textbf{Convention:} \textit{All over this paper, we assume that provability models are classical.}

\end{definition}

\begin{remark}
As we see, all provability models are assumed to be  \textit{classical}, 
i.e.~the theories attached to the worlds are assumed to be classical. Nevertheless, 
one may consider them to be not being necessarily classical. For example, 
in \cite{PLHA}, the first author considers attached theories as intuitionistic modal theories. 
Also \cite{Mayaux} considers them as a mixture of classical and intuitionistic theories. 
\end{remark}

When proving soundness results it will be common to implicitly use the following result, 
which allows us to skip half of the proof.
\begin{lemma}
Let $A$ be a purely modal formula and $\pcal$ be a provability model.
If $\pcal \pmodels A$ then for any $w \in W^\R$, $\vdash_w  A$.
\end{lemma}
\begin{proof}
    Let $w \in W^\R$, since $\pcal \pmodels A$ we have that $\pcal, w \pmodels^+ A$.
    Using that $A$ is purely modal and by modal completeness we obtain the desired 
    $\vdash_w A$.
\end{proof}

\begin{remark}\label{lcali-lcalb-dist}
Notice that the above notions
rely on the specific language of modal logic which we are working with. 
For example,  the logical properties like being classical, depends on the language in the following sense. 
Classicality of provability models for the language $\lcalb$ 
implies that $A\to(B\to A)$ should belong 
to the set of axioms of $\lgc w$ for every $w\in W^\R$ and $A,B\in\lcalb$, 
while for  the language 
$\lcali$ (which will be considered in \cref{sec3}), it implies that  
$A\to(B\to A)$ should belong to the set of axioms of $\lgc w$ 
for every $w\in W^\R$ and $A,B\in\lcali$.
Another example is the notion of modal completeness. It relies on the notion of 
purely modal, which varies in different modal languages.
Having these differences in mind, 
we keep using these notions in future sections for the other 
modal languages as well. 
Nevertheless, 
to avoid confusion in our future uses, 
we make it clear that for which language we are considering 
a provability model: $\lcalb$, $\lcali$ or $\lcalo$. 
More explicitly, we always consider provability models 
in \cref{sec2,sec2-2} for $\lcalb$, in 
\cref{sec3} we always consider them for $\lcali$ and in 
\cref{sec-GLP} we always consider them for $\lcalo$.
\end{remark}

Before  we continue with the soundness and completeness of $\sfk$, 
let us have some general correspondence
between provability models and Kripke models. This shows that the class of 
provability models, \textit{includes} the class of all Kripke models.

\begin{lemma}\label{lem1}
For every Kripke model $\kcal=(W,\R,\V)$, there is a provability model 
$\pcal=(W, \R, \{\lgc w\}_{w\in W},\V)$ which is equivalent to 
$\kcal$ in the following sense:
$$\forall A\, \forall w\in W \ (\kcal,w\models A \ \Leftrightarrow \ \pcal,w\pmodels A).$$
Moreover, if $\kcal$ is transitive, then we may pick 
$\pcal$ having necessitation.
\end{lemma}
\begin{proof}
For $\kcal$ without transitivity,
let $\axioms{\lgc w}:=\{A: \kcal,w\models A\}$ and $\rules{\lgc w}$ 
to include only the modus ponens. 
It is clear that $\axioms{\lgc w}$ has all classical tautologies. Hence 
$\pcal$ is classical.
%Also by definition, $\pcal$ is modally complete, and hence a provability model indeed. 
Furthermore, observe that
$\axioms{\lgc w}$ is already closed under modus ponens. Hence 
$\vdash_w A$ is equivalent to $\kcal,w\models A$. 

By induction on the complexity of formulas we show that this provability model 
is equivalent to $\kcal$, i.e., that 
$\forall w \in W \ (\kcal, w \models A \Leftrightarrow \pcal,w\pmodels A)$.
The only interesting case is when $A = \Box B$.
Assume $\kcal, w \models \Box B$, then  for any 
$v \sqsupset w$ we have that $\kcal, v \models B$.
By definition this means $\forall v \sqsupset w \ (\, \vdash_v B)$, 
i.e., $\pcal, w \pmodels \Box B$.
For the opposite direction, assume that $\pcal, w \pmodels \Box B$, 
so $\forall v \sqsupset w \ (\,\vdash_v B)$.
By definition of $\lgc v$ we obtain that $\forall v \sqsupset w (\kcal, v \models B)$, 
i.e., $\kcal, w \models \Box B$.
Finally, we can conclude modal completeness of $\pcal$ easily by using the equivalence 
of $\kcal$ and $\pcal$.

For $\kcal$ with transitivity, 
let $\axioms{\lgc w}:= \{ A: \kcal,w\models \Boxdot A\}$ and 
$\rules{{\lgc w}}$ to include two inference rules: modus ponens and necessitation.  
By definition, it is easy to check that $\axioms{\lgc w}$ 
contains all propositional tautologies and it is already closed under modus ponens.
Additionally, thanks to transitivity, $\axioms{\lgc w}$  is also closed 
under necessitation as well. This means that $\vdash_w A$ is 
equivalent to $\kcal,w\models \Boxdot A$.
By induction on the complexity of formulas we show that this provability premodel 
is equivalent to $\kcal$, i.e., that 
$\forall w \in W\  (\kcal, w \models A \Leftrightarrow \pcal,w\pmodels A)$.
The only interesting case is when $A = \Box B$.
Assume $\kcal, w \models \Box B$, 
then (by transitivity) for any $v \sqsupset w$ we have that 
$\kcal, v \models \Boxdot B$.
By definition this means $\forall v \sqsupset w\ (\lgc v \vdash B)$, i.e., 
$\pcal, w \pmodels \Box B$.
For the opposite direction assume that $\pcal, w \pmodels \Box B$, so 
\(\forall v \sqsupset w\  (\lgc v \vdash B)$.
By our earlier observation we obtain that 
\(\forall v \sqsupset w\ (\kcal, v \models \Boxdot B)$, 
so in particular $\forall v \sqsupset w\ (\kcal, v \models B)$, 
i.e.~$\kcal, w \models \Box B$.
Finally, we can conclude modal completeness of $\pcal$ 
easily by using the equivalence of $\kcal$ and $\pcal$.
\end{proof}

\subsection{Provability models for $\sfk$}
In this subsection, we show that $\sfk$ is sound and complete for provability models. 

For the proof of \Cref{soundness-k}, we use the following alternative axiomatization of 
the modal logic $\sfk$, in which we don't have the necessitation rule and instead 
we have $\Box^n$ of every axiom for every $n\in\nat$.
Thus $\sfk$ has only modus ponens as its inference rule and the following list of axioms
for every $n\in\nat$ (take $\Box^0 A:=A$.):
\begin{itemize}
\item $\Box^n A$ for 
every classical tautology $A$.
\item $\Box^n\big(\Box(A\to B)\to (\Box A\to \Box B)\big)$.
\end{itemize}

It is well-known in the literature that $\sfk$ is sound and complete for 
finite Kripke models.
\begin{theorem}\label{soundness-k}
$\sfk$ is sound for provability models.
\end{theorem}
\begin{proof}
Let $\pcal=(W,\R,\{\lgc w\}_{w\in W},\V)$ be a provability  model.
It is enough to show that all axioms of $\sfk$, are true in $\pcal$. 
\begin{itemize}
\item $A=\Box^{n} B$ and $B$ is classical tautology: We use induction on $n$
to prove this. The case $n=0$ is immediate. So assume that we already have 
the validity of $\Box^n B$ for every classical tautology $B$ at any $u\in W$, seeking to 
show $\pcal,w\pmodels \Box^{n+1} B$. 
By definition it is enough to show $\vdash_u \Box^n B$ for every $u\sqsupset w$.
If $n=0$, then since 
$\lgc u$ includes all classical tautologies we have desired result. If also 
$n> 0$, then by induction hypothesis we have $\pcal,u\pmodelsp \Box^n B$ and since 
$\Box^n B$ is purely modal, by modal completeness we get $\vdash_u \Box^n B$, as desired. 
\item $A=\Box^n\big(\Box(B\to C)\to(\Box B\to \Box C)\big)$. Again we use induction on $n$.
The case $n=0$ is sound because each $\lgc u$ includes modus ponens.
So assume as induction hypothesis that we have $\pcal,w\pmodels \Box^n\big(\Box(B\to C)\to(\Box B\to \Box C)\big)$ for every $w\in W$, seeking to show 
$\pcal,w\pmodels \Box^{n+1}\big(\Box(B\to C)\to(\Box B\to \Box C)\big)$. 
By definition, it is enough to show $\vdash_u \Box^n\big(\Box(B\to C)\to(\Box B\to \Box C)\big)$, which holds by induction hypothesis and 
the fact that $\Box^n\big(\Box(B\to C)\to(\Box B\to \Box C)\big)$ 
is purely modal (even when $n=0$).\qedhere
\end{itemize}
\end{proof}

\begin{theorem}\label{completeness-k}
$\sfk$ is complete for finite provability models.
\end{theorem}
\begin{proof}
It is a direct consequence of \Cref{lem1} and 
the completeness of $\sfk$ for finite Kripke models. 
\end{proof}

\subsection{Provability models of $\kfour$}

The necessitation rule, is defined as \Ax{$A$}\UI{$\Box A$}\DP.
For provability models of $\kfour$ we have the flexibility to consider either 
of following properties on the provability models:
\begin{itemize}
\item The frame property,  transitivity, or 
\item the logical property necessitation, i.e., that for any $w \in W^\R$ if $\lgc w \vdash A$ then $\lgc w \vdash \Box A$.
\end{itemize}
This means that, the soundness and completeness of $\kfour$ holds both for 
transitive provability models and provability models with logical property of 
closure under necessitation.

Again we consider an equivalent axiomatization of $\kfour$ without necessitation. 
In other words, $\kfour$ includes only modus ponens as its inference rules. 
As axioms of $\kfour$ we have the following:
\begin{itemize}
\item All above-mentioned axioms of $\sfk$\footnote{Actually, in the presence of the transitivity axiom, it is enough to consider 
the axioms of $\sfk$ with $n\leq 1$.}.
\item $\Box A\to\Box\Box A$.
\item $\Box (\Box A\to\Box\Box A)$.
\end{itemize}

\begin{theorem}\label{soundness-k4}
$\kfour $ is sound for provability models with necessitation.
Also it is sound for transitive provability models.
\end{theorem}
\begin{proof}
Let  $\pcal=(W,\R,\{\lgc w\}_{w\in W},\V)$ be a provability model with necessitation.
By \Cref{soundness-k} it is enough only to show the validity of following
axioms of $\kfour$:
\begin{itemize}
\item $\Box A\to \Box \Box A$: Since $\lgc w$ has necessitation.
\item $\Box (\Box A\to \Box \Box A)$: by modal completeness and above item.
\end{itemize}
Next assume that $\pcal$ is a transitive provability model, so ${\R^+} = {\R}$. 
Again, it is enough to show that all axioms are valid in $\pcal$. 
We only treat the following case:
\begin{itemize}
\item $\Box A\to \Box\Box A$: Let $\pcal,w\pmodels \Box A$ and $u\sqsupset w$.
We need to show $\lgc u \vdash \Box A$. For this, it suffices to show 
$\pcal, u \pmodels^+ \Box A$, and use modal completeness.
We pick $w$ to show that $\pcal, u \pmodels^+ \Box A$, i.e., we need to prove that
$\forall v \sqsupset^+ w\  (\pcal, v \pmodels \Box A)$.
Pick $v \sqsupset^+ w$ and $v' \sqsupset v$ arbitrary, then $w \R^+ v'$ and by transitivity we have $w \R v'$.
Then $\pcal, w \pmodels \Box A$ implies $ \vdash_{v'} A$, so $\pcal, v \pmodels \Box A$ as desired.
\qedhere
\end{itemize}
\end{proof}

\begin{theorem}\label{completeness-k4}
$\kfour$ is complete for finite transitive provability models 
with necessitation.
\end{theorem}
\begin{proof}
It is a direct consequence of \Cref{lem1} and 
the completeness of $\kfour$ for finite transitive Kripke models. 
\end{proof}

\begin{remark}
    Note that thanks to Theorem~\ref{soundness-k4} and Theorem~\ref{completeness-k4} we have that $\four$ is sound and complete with respect to the class of transitive provability models and to the class of provability models with necessitation.
    In other words, only one of transitivity or necessitation is needed.
\end{remark}

\subsection{Provability models of $\sfour$}

As we will show in this section, the modal logic 
$\sfour$ is sound and complete for provability models 
of $\kfour $ with following additional properties:
\begin{itemize}
\item The frame property reflexivity.
\item The logical property, \textit{local completeness}: 
$\pcal,w\pmodelsp A$ implies $\vdash_w A $ for every $w\in W^\R$ and 
$A\in\lcalb$.\footnote{Note the difference with modal completeness. 
In local completeness we do not require the formula to be purely modal.}
\item The logical property of \textit{local soundness}:  $\vdash_w A $
implies $\pcal,w\pmodels A$, for every $A$ and $w\in W$.
\end{itemize}

Again we consider an equivalent axiomatization of $\sfour$ without necessitation. 
In other words, $\sfour$ defined as $\kfour$ plus the following axioms:
\begin{itemize}
\item $\Box A\to A$.
\item $\Box (\Box A\to A)$.
\end{itemize}

\begin{theorem}\label{soundness-s4}
$\sfour $ is sound for reflexive 
provability models with necessitation and 
local soundness and local completeness.
 Moreover, the necessitation can be replaced by transitivity. 
\end{theorem}
\begin{proof}
Let  $\pcal=(W,\R,\{\lgc w\}_{w\in W},\V)$ be a reflexive provability model with necessitation, local soundness and local completeness.
By \Cref{soundness-k4} it is enough only to show the validity of following
axioms of $\sfour$:
\begin{itemize}
\item $\Box A\to A$: 
Immediate by reflexivity and local soundness of 
$\lgc w$.
\item $\Box (\Box A\to  A)$: by local completeness and above item.\footnote{Note that modal completeness does not suffice for this, as $\Box A \to A$ is not purely modal.}
\end{itemize}
The proof of soundness when we have transitivity instead of necessitation, is in the same way as for $\kfour$.
\end{proof}

\begin{theorem}\label{completeness-s4}
$\sfour$ is complete for finite transitive reflexive provability models 
with necessitation, local soundness and local completeness.
\end{theorem}
\begin{proof}
Let $\sfour \not\vdash A$, by the completeness of $\sfour$ for finite reflexive 
and transitive Kripke models there is such a Kripke model $\kcal$ with 
$\kcal, w \not\models A$. 
Consider the provability model $\pcal$ obtained from $\kcal$ as given by \Cref{lem1}.
By construction we know it is finite, transitive, reflexive, and  with necessitation. Also
by \Cref{lem1} we know that $\pcal, w \not\pmodels A$. So all left to notice is 
that $\pcal$ is indeed locally complete and sound.

Assume $\pcal, w \pmodels^+ A$, then $\pcal, w \pmodels \Boxdot A$ so 
$\kcal, w \models \Boxdot A$ and then $\vdash_w A$, so $\pcal$ is locally complete.
Now, assume that $\vdash_w A$. By definition $\kcal,w \models \Boxdot A$ so $\kcal, w \models A$ and then $\pcal, w \pmodels A$, so $\pcal$ is locally sound.
\end{proof}

A closer look at the provability models of $\sfour$, reveals that we actually can not 
get far from the standard Kripke models of $\sfour$ in this case. 
The reason  is for the 
combination of local soundness and local completeness. This combination, actually
enables us to recover a Kripke model of $\sfour$, out of a given provability model. This is the converse of \Cref{lem1}, which shows that the provability models of $\sfour$
are actually equivalent to Kripke models.

\begin{lemma}\label{lem2}
For every transitive provability model $\pcal=(W,\R,\{\lgc w\}_{w\in W},\V)$  
  with local soundness and completenesses, the  Kripke model 
$\kcal_\pcal:=(W,\R,\V)$ is equivalent with $\pcal$. More precisely,
\begin{center}
$\pcal,w\pmodels A$ iff $\kcal_\pcal,w\models A$ for any $w\in W$ and $A\in\lcalb$.
\end{center}
\end{lemma}
\begin{proof}
We use induction on $A$. All cases are trivial, except $A=\Box B$. 
First let $\kcal_\pcal,w\nmodels\Box B$ and take $u\sqsupset w$ such that 
$\kcal_\pcal,u\nmodels B$.
Then by induction hypothesis, we have $\pcal,u\npmodels B$. This by local soundness implies
$B\nin\lgc u$, and hence $\pcal,w\npmodels \Box B$.

For the other way around, let $\kcal_\pcal,w\models \Box B$ and $u\sqsupset w$.
We need to show $\vdash_u B$. Since $\kcal_\pcal,w\models \Box B$, 
by transitivity, we get $\kcal_\pcal,u\modelsp B$. Hence induction hypothesis 
implies $\pcal,u\pmodelsp B$, and thus by local modal completeness we have $\vdash_u B$.
\end{proof}

\subsection{Provability models of $\GL$}\label{sec-gl}
In this section we show that $\GL$ is sound and complete for converse well-founded 
provability models with necessitation.  Similar to $\kfour$, we have some flexibilities:
\begin{itemize}
\item Converse well-foundedness $\Leftrightarrow$ L\"ob's Rule.
\item Transitivity $\Leftrightarrow$ Necessitation.
\end{itemize}
By L\"ob's rule we mean 
\Ax{$\Box A\to A$}
\UI{$A$}
\DP.

Again we consider an equivalent axiomatization of $\GL$ without necessitation. 
In other words, $\GL$  is defined as $\kfour$ (as defined earlier in this paper)
plus the following axioms:
\begin{itemize}
\item L\"ob's axiom: $\Box (\Box A\to A)\to \Box A$.
\item  The box of L\"ob's axiom: $\Box B$ for $B=\Box (\Box A\to A)\to \Box A$.
\end{itemize}
\begin{theorem}\label{soundness-gl}
$\GL $ is sound for   provability models having  necessitation and L\"ob's rule. 
Moreover, the soundness still holds if we do either of replacements (or both):
\begin{itemize}
\item Necessitation replaced by transitivity of frame,
\item L\"ob's rule replaced by converse well-foundedness.
\end{itemize}
\end{theorem}
\begin{proof}
Let  $\pcal=(W,\R,\{\lgc w\}_{w\in W},\V)$ be a provability model with either 
necessitation or transitivity. Also assume that $\pcal$ is either converse 
well-founded or having L\"ob's rule as its logical property.
By \Cref{soundness-k4} it is enough only to prove  the soundness of L\"ob's axiom and 
its box. The soundness of its box holds by the property of modal completeness. 
So we only show the soundness  of L\"ob's axiom in two cases:
\begin{itemize}
\item $\pcal$ is converse well-founded: By induction on $w\in W$ we show that 
L\"ob's axiom is valid at $w$. So as induction hypothesis we assume that  
L\"ob's axiom is valid at any $u\sqsupset w$. 
In particular, induction hypothesis implies that 
$\pcal,u\pmodelsp \Box(\Box A\to A)\to\Box A$ for every $u\sqsupset w$. Hence by modal completeness,
$\lgc u$ implies the L\"ob's axiom
for every $u\sqsupset w$.
Let $\pcal,w\pmodels \Box(\Box A\to A)$.
We must show that $\pcal,w\pmodels \Box A$. So let $u\sqsupset w$. Then we have 
$ \vdash_u \Box A\to A$. By necessitation, this implies 
$\vdash_u \Box (\Box A\to A)$. Then induction hypothesis and 
modus ponens implies $\vdash_u \Box A$. This together with $\vdash_u \Box A\to A$
implies $\vdash_u A$. Thus $\pcal,w\pmodels \Box A$, as desired. 
\item $\pcal$ has the logical property of L\"ob's rule: 
Let $\pcal,w\pmodels \Box(\Box A\to A)$.
We must show that $\pcal,w\pmodels \Box A$. So let $u\sqsupset w$. Then we have 
$\vdash_u \Box A\to A$. By  L\"ob's rule this implies  $\vdash_u   A$.   
Thus $\pcal,w\pmodels \Box A$, as desired.
\qedhere
\end{itemize}
\end{proof}

Thanks to the completeness of $\GL$ for finite irreflexive transitive models we have the following corollary.
\begin{corollary}\label{finfalsGL}
$\GL\vdash A$ iff $\GL\vdash \Box^k\bot\to A$ for all 
$k\in\nat$.
\end{corollary}

Let us define $\GL_n :=\GL +\Box^n\bot$, i.e., adding $\Box^n \bot$ to the axioms of $\GL$ and closing them under modus ponens.
It is easy to show that $\GL_n$ is closed under necessitation.
\begin{lemma}\label{gl-local-tabularity}
For every $n$, $\GL_n$ is locally tabular, i.e.~there are finitely many 
formulas over the finite language modulo $\GL_n$-provable equivalence. 
Moreover, One may calculate a finite set $X_n$ which includes a  $\GL_n$-equivalent of 
every formula.
\end{lemma}
\begin{proof}[Sketch of the proof.]
Given that Kripke models of $\GL_n$ are the finite irreflexive transitive models with height $n$, one may easily observe that every formula is equivalent to a formula in which the maximum number of nesting 
$\Box$'s are $n$. Then one may easily construct the desired finite set by induction on $n$.
\end{proof}

\begin{theorem}\label{completeness-gl}
$\GL$ is complete for finite transitive irreflexive provability models 
having the following logical properties: necessitation and L\"ob's rule.
\end{theorem}
\begin{proof}
It is a direct consequence of \Cref{lem1} and 
the completeness of $\GL$ for finite irreflexive transitive Kripke models. 
\end{proof}

\subsection{Strong completeness of $\GL$ for provability models}
It is known that $\GL$ is not strongly complete for its standard Kripke models \cite[Theorem 4.43]{Blackburn}.
Nevertheless, it is strongly complete for its topological semantics \cite{Esa81,aguilera2017strong}. 
Here we show that $\GL$ is also strongly complete for provability models.\footnote{We remember 
that given a class of models $\pfrak$ and a theory $\lgc{}$ (where $\lgc{}$ 
represents an axiomatization with only the rule of modus ponens) we say that $\lgc{}$ 
is strongly complete for $\pfrak$ if $\Gamma \pmodels_{\pfrak} A$ implies 
$\Gamma \vdash_\sft A$.}

%The results of this subsection are also interesting for the following reason. 
%Up to this point, we have considered provability models for various modal theories. 
%However, as indicated by the proofs of completeness theorems, 
%these models are not genuinely 
%\textit{new} tools independent of standard Kripke models. 
%So far, the only method to construct a provability model has been to start with a 
%standard Kripke model, assign to each node a set of valid formulas 
%(as we described in the proof of \Cref{lem1}), and thereby obtain 
%an equivalent provability model. In contrast, the following theorem demonstrates the 
%existence  of a provability model that is genuinely \textit{independent} of standard 
%Kripke models. By independence, we mean that such a model is not merely an equivalent 
%provability model of some standard Kripke model, as defined in the proof of 
%\Cref{lem1}. This distinction arises from the fact that $\GL$ is not strongly 
%complete with respect to standard Kripke models of $\GL$ 
%(see \cite[Theorem 4.43]{Blackburn}).

\begin{theorem}\label{strong-comp}
$\GL$ is strongly complete for the provability 
models with necessitation and L\"ob's rule.
\end{theorem}
\begin{proof}
Let $\Gamma$ be a set of formulas such that 
$\Gamma\nvdash_\GL A$. Recall from \cref{Fin-seq}, the notations on finite sequences.
Define 
\begin{itemize}
\item $W_0$ as the set of all maximal consistent  
sets of formulas including $\GL$. 
%By maximal consistent, we mean a set $\Gamma'$ such that 
%$\Gamma'\nvdash\bot$ and $\Gamma'\vdash\neg B$ for every $B\nin\Gamma'$.
\item $w\prec u$ for every $w,u\in W_0$ iff
for all   
$\Box B\in w$ we have $B\in u$.
\item $W:=\{\langle w_0,\ldots,w_n\rangle: 
	w_0\prec\ldots\prec w_n\}$.
\item $\sigma\R \tau$ iff $\sigma,\tau\in W$
and there is some $w\in W$	such that $\tau=\sigma*\langle w\rangle$.
\item $\sigma\V p$ iff $p\in \fel \sigma $.
\item For $\sigma:=\iseg\tau$
define $\lgcp\tau$ to include axioms
$\axioms{\lgcp \tau}:=\{B\in\lcalb: \forall\,u\succ \fel\sigma\ (B\in u)\}$ and 
without any inference rules.  This means that ${\lgcp \tau}\vdash B$ is equivalent to 
$B\in\axioms{\lgcp\tau}$.
Then define $\lgc\tau$ to be $\lgcp\tau$ augmented with the inference rules: 
modus ponens, necessitation and L\"ob's rule.
%As we will see later in claims 2 and 3, the axiom set
%of $\lgc \tau$ is already closed under the modus ponens, necessitation and L\"ob's rule.
%Hence without any change to the validity of formulas in the pre-provability model, we may 
%replace $\lgcp\tau$ with $\lgc\tau$.
\item Finally define 
$\pcal':=(W,\R, \{\lgcp\sigma\}_{\sigma\in W^{\R}},\V)$
and $\pcal:=(W,\R, \{\lgc\sigma\}_{\sigma\in W^{\R}},\V)$.
\end{itemize}
Obviously both $\pcal'$ and $\pcal$ are pre-models. 
\\[2mm]
\textit{Claim 1.} $\pcal',\sigma\pmodels B$ iff $B\in\fel\sigma$ for every $B\in\lcalb$ and $\sigma\in W$.
\\[2mm]
\textit{Proof of the claim 1.} We prove this by induction on $B$. All cases are obvious except for 
$B=\Box C$. First let 
$\Box C\in \fel\sigma$. We want to show 
$\pcal',\sigma\pmodels \Box C$.  Let $\tau=\sigma*\langle w\rangle\in W$. 
It is enough to show that $ C\in\axioms{\lgcp\tau}$.  Since $\Box C\in \fel\sigma$, 
then for any $u\succ\fel\sigma$ we have 
$C\in u$ and hence by definition we have $C\in\axioms{\lgcp\tau}$, as desired.
For the other way around, assume that $\Box C\nin\fel \sigma=w$
seeking some $\tau\sqsupset \sigma$ such that $\nvdash_{\lgcp \tau} C$.
By definition of $\lgcp \tau$ and $\sigma\R\tau$, 
it is enough to find some $u\succ w$ such that $C\nin u$.
Define $u_0:=\{E\in\lcalb: \Box E\in w\}$.
Notice that $\GL\subseteq u_0$ and $u_0\nvdash C$. 
Thus there is a maximal set $u\supseteq u_0$ such that $ C\nin u$. Given that 
$u\supseteq u_0$, we get  $w\prec u$, as desired.
\\[2mm]
\textit{Claim 2.} $\lgcp\sigma$ is closed under modus ponens, 
includes axioms of classical logic and has modal completeness.
\\[2mm]
\textit{Proof of the claim 2.}
First notice that since each $w\in W_0$, by maximality, includes classical logic, 
then classical axioms are already included in $\axioms{\lgcp \sigma}$ for any 
$\sigma\in W^\R$.
Moreover, since each $w\in W_0$ is closed under modus ponens (again by maximality), 
then $\lgcp \sigma$ is also closed under modus ponens.
We only need to show modal completeness.  So assume that $E$ is purely modal such that 
$\pcal',\sigma\pmodelsp E$ for some $\R$-accessible 
 $\sigma\in W$. Then by claim 1 above, we have 
 $E $ belongs to the set of axioms of $\lgcp\sigma$, as desired.
\\[2mm]
\textit{Claim 3.} $\lgcp\sigma$ is closed under necessitation and 
L\"ob's rule.
\\[2mm] 
\textit{Proof of the claim 3.}  
For the necessitation, let $B\in\axioms{\lgcp\sigma}$ for 
$\sigma=\tau*\langle w\rangle\in W$. 
This implies by definition that for any $u\succ \fel\tau$ we have 
$B\in u$. It is enough to show that  for any $u\succ \fel\tau$ 
also we have $\Box B\in u$. Let $u\succ \fel\tau$ and take 
$\eta:=\tau*\langle u\rangle$  and  $\pi\sqsupset \eta$. 
Then obviously $B\in\axioms{\lgcp \pi}$ and hence by definition, 
$\pcal',\eta\pmodels \Box B$. This means that 
for every $\eta\sqsupset \tau$ we have $\pcal',\eta\pmodels \Box B$. This together with 
claim 1 implies that $\Box B\in\axioms{\lgcp \sigma}$, as desired.

For the L\"ob's rule, let $\Box B\to B\in\axioms{\lgcp\sigma}$. 
By necessitation (as proved just above) 
we have $\Box(\Box B\to B)\in \axioms{\lgcp\sigma} $. Then since 
L\"ob's axiom holds in all $w\in W_0$, we have it also 
in $\axioms{\lgcp \sigma}$. 
Thus by modus ponens (as proved to hold in claim 2) we have 
$\Box B\in\axioms{\lgcp \sigma}$. 
By another use of modus ponens, then we have $B\in\axioms{\lgcp \sigma}$, as desired.
\\[2mm]
\textit{Claim 4.} There is a $\sigma\in W$ such that 
$\pcal,\sigma\npmodels A$
while $\pcal,\sigma\pmodels B$ for every $B\in \Gamma$. 
\\[2mm]
\textit{Proof of the claim 4.}
First notice that by previous claims, replacing $\lgc\tau$ for $\lgcp \tau$, does not affect the evaluation of formulas in the provability model. Thus $\pcal'$ 
and $\pcal$ are sharing the same validity for formulas. 
%Then, together with previous claims, it finishes the proof of this Theorem.
So to prove this claim, we reason as follows. 
Since $\Gamma,\GL\nvdash A$, there is a maximal consistent set $w\supseteq \Gamma\cup\GL$ such that 
$A\nin w$. Then claim 1 implies $\pcal,\langle w\rangle \npmodels A$ and $\pcal,\langle w\rangle \pmodels \Gamma$, as desired.
\end{proof}

\section{Construction of Provability Models}\label{sec2-2}
%Although in \Cref{strong-comp} we gave a
%canonical provability model for $\GL$ for which it is strongly complete, it is still not 
%satisfying for a picky reader, given that this model is not defined in a constructive way. 
%Even the decidability of $\pcal,w\pmodels B$ is not quarantined. 
%
%In this section, we show that it is possible to go beyond and define 
%independently a provability model in the case of $\GL$ (\Cref{bases-gl}). 
%More precisely, for every converse well-founded tree frame, and every assignment of sets
%of formulas to the nodes, one may expand the sets such that the result becomes 
%a provability model with necessitation. 
%
%Another defect in provability models is the following. Validity of formulas in 
%provability models, might not be decidable. 
%This may have arisen due to the undecidability of the theories associated with the nodes.
%However, as we will see in \Cref{bases-gl}, the constructed provability model, satisfies 
%decidability, in case of finiteness of the sets being assigned to the worlds.
%
In previous section, we have seen   provability models for 
various modal logics. 
Nevertheless, as can be seen through the proofs of completeness theorems, 
they are not looking as 
independent tools from standard Kripke models. 
The only way to build a provability model, up to now, is by taking a standard 
Kripke model and attaching a set of valid formulas to the node and finally 
making an equivalent provability model. More precisely, for all constructed Provability 
models $\pcal:=(W,\R, \{\lgc\sigma\}_{\sigma\in W^{\R}},\V)$, the standard part of the 
model $\kcal_\pcal:=(W,\R, \V)$, is equivalent to $\pcal$ in the following sense: 
$$\pcal,w\pmodels A \quad\Leftrightarrow \quad \kcal_\pcal,w\pmodels A.$$

In this section, we show that it is possible to go beyond and define 
independently a provability model in the case of $\GL$ (\Cref{bases-gl}). 
More precisely, for every converse well-founded tree frame, and every assignment of sets
of formulas to the nodes, one may expand the sets as theories such that the result becomes 
a provability model with necessitation\footnote{Of course one may easily define an 
independent provability model, by attaching inconsistent theories to the nodes. 
However, importance of the content of this section is that we do it in a systematic way, 
for arbitrary assignments to the nodes.}. 

Another defect in provability models is the following: 
Validity of formulas in 
provability models, might not be decidable 
(as it was the case for instance in the canonical provability model constructed for 
$\GL$ in \Cref{strong-comp}, for which, we had strong completeness). 
This might have arisen due to the undecidability of the theories associated with the nodes.
However, as we will see in \Cref{bases-gl}, the constructed provability model, satisfies 
decidability, in case of finiteness of the sets being assigned to the worlds.

\begin{definition}
    Let $A$ be a modal formula whose atomic subformulas 
     are $\vec{p} = p_0,\ldots,p_{m-1}$.
    A $\Box$-free formula $B(\vec{p},q_0,\ldots,q_{n-1})$ is said to be a
    \emph{propositional skeleton of $A$} if the atomics $q_0,\ldots,q_{n-1}$ 
    are distinct from $\vec{p}$ and there are unique and distinct 
    modal formulas $C_0,\ldots,C_{n-1}$ such that
    $A = B(\vec{p}, \Box C_0,\ldots,\Box C_{n-1})$.
    Then, we define the purely modal uniform pre-interpolant of $A$ denoted as $A^*$ as
    \[
    A^*:=
    \bigwedge_{\vec{b} \in \{\top, \bot\}^m} B(\vec{b}, \Box C_0, \ldots, \Box C_{n-1}). 
    \]
\end{definition}
We note the two following facts:
\begin{enumerate}
    \item Every modal formula has a propositional skeleton, although it is not unique (nevertheless, this nonuniqueness is only due to renaming of the variables $q_0,\ldots,q_{n-1}$).
    \item Even if the propositional skeleton of a formula is not unique, 
    the purely modal uniform interpolant is unique (once we fix an order on the conjuncts).
\end{enumerate}
Before proving the properties of purely modal uniform pre-interpolation we need the following simple lemma\footnote{We remember that $\vdash$ means provability in classical 
logic, but formulated in the modal language.}.
%Notice that $\vdash$ is used as derivability in classical logic without any modal axioms.
\begin{lemma}\label{lm:boxes-are-variables}
    Let $A(p_0,\ldots,p_{m-1})$ be a $\Box$-free formula and $B_0,\ldots,B_{m-1}$ be modal formulas.
    Then $\vdash A(\Box B_0, \ldots, \Box B_{m-1})$ implies $\vdash A(p_0,\ldots,p_{m-1})$.
\end{lemma}
\begin{proof}
    By induction in the proof of $A(\Box B_0, \ldots, \Box B_{m-1})$, 
    since with classical propositional tautologies and modus ponens $\Box$-formulas 
    only act as atomic formulas.
\end{proof}

The properties of modal uniform pre-interpolant that we are going to need (and the justification for this name) are collected in the following lemma.

\begin{lemma}
    Let $A$ be a modal formula.
    We have the following:
    \begin{enumerate}
        \item $A^*$ is purely modal.
        \item $\vdash A^* \to A$.
        \item For any purely modal formula $B$, we have that $\vdash B \to A$ implies 
        $\vdash B \to A^*$.
    \end{enumerate}
\end{lemma}
\begin{proof}
The first property is trivial by definition.
Let us show the second property. 
Let $A$ be a formula whose propositional variables are $p_0,\ldots,p_{m-1}$ and $C(\vec{p},q_0,\ldots,q_{n-1})$ be a propositional skeleton of $A$, such that $A = C(\vec{p}, \Box B_0,\ldots, \Box B_{m-1})$.
Note that $\bigwedge_{\vec{b} \in \{\bot,\top\}^m} C(\vec{b},q_0,\ldots,q_{m-1})$ is the uniform pre-interpolant in classical propositional logic of $C(\vec{p}, q_0,\ldots,q_{m-1})$ for the language $\{q_0,\ldots,q_{m-1}\}$.
In particular, this means that $\vdash\bigwedge_{\vec{b} \in \{\bot,\top\}^m} C(\vec{b},q_0,\ldots,q_{m-1}) \to C(\vec{p},q_0,\ldots,q_{m-1})$.
Since classical propositional logic is closed under substitution and the variables in $\vec{p}$ and $q_0,\ldots,q_{m-1}$ are different, we obtain that $\vdash \bigwedge_{\vec{b} \in \{\bot,\top\}^m} C(\vec{b},\Box B_0,\ldots,\Box B_{m-1}) \to C(\vec{p},\Box B_0,\ldots,\Box B_{m-1})$, i.e., $\vdash A^* \to A$.

Finally, we show the third property.
Assume that $\vdash B \to A$ where $B$ is purely modal.
Then, we can get a propositional skeleton $C_B(q_0,\ldots,q_{m-1})$ of $B$  and a propositional skeleton $C_A(\vec{p}, q_0,\ldots,q_{m-1})$ of $A$ such that
\[
B = C_B(\Box B_0,\ldots,\Box B_{m-1}) \qquad \text{ and } \qquad A = C_A(\vec{p},\Box B_0,\ldots,\Box B_{m-1}).
\]
In other words, we search the propositional skeleton such that the extra $q$-variables are shared between $C_B$ and $C_A$ in case they need to be substituted by the same formula.
Then the assumption says $\vdash C_B(\Box B_0,\ldots,\Box B_{m-1}) \to  C_A(\vec{p},\Box B_0,\ldots,\Box B_{m-1})$, where $C_B(q_0,\ldots,q_{m-1}) \to C_A(\vec{p}, q_1,\ldots,q_{m-1})$ is a $\Box$-free formula.
By Lemma~\ref{lm:boxes-are-variables} we have that $\vdash C_B(q_0,\ldots,q_{m-1}) \to C_A(\vec{p}, q_1,\ldots,q_{m-1})$.
Since $\bigwedge_{\vec{b} \in \{\bot,\top\}^m} C_A(\vec{b},q_0,\ldots,q_{m-1})$ is the uniform pre-interpolant in classical propositional logic of $C_A(\vec{p}, q_0,\ldots,q_{m-1})$ for the vocabulary $\{q_0,\ldots,q_{m-1}\}$ and the vocabulary of $C_B$ is contained in $\{q_0,\ldots,q_{m-1}\}$, we obtain that $\vdash C_B(q_0,\ldots,q_{m-1}) \to \bigwedge_{\vec{b} \in \{\bot,\top\}^m} C_A(\vec{b},q_0,\ldots,q_{m-1})$.
Since classical propositonal logic is closed under substitution, we obtain $\vdash C_B(\Box B_0,\ldots,\Box B_{m-1}) \to \bigwedge_{\vec{b} \in \{\bot,\top\}^m} C_A(\vec{b},\Box B_0,\ldots,\Box B_{m-1})$, i.e., the desired $\vdash B \to A^*$.
\end{proof}

Let $\pcal$ and $\pcal'$ be two provability pre-models. We say that $\pcal$
is a \textit{submodel} of $\pcal'$, notation 
$\pcal\subseteq\pcal'$, if 
$\pcal$ and $\pcal'$  are sharing the same frames, the same valuations on variables, 
while for every $w\in W$, the logic
$\lgc w ' $ assigned to $w$ in $\pcal'$ includes the theory $\lgc w$ being assigned to 
$w$ in $\pcal$\footnote{A theory $\sft'$ includes $\sft$ if every theorem of $\sft$ 
is also a theorem of $\sft'$.
In other words, ${\sft}\vdash A$ implies 
${\sft'}\vdash A$. 
%both the set of axioms and rules 
%of inferences of $\sft'$ are including the axioms and rules of inferences in $\sft$, 
%i.e.~$\axioms{\sft'}\supseteq\axioms{\sft}$ and $\rules{\sft'}\supseteq\rules{\sft}$
.}.
Notice that   $\pcal\subseteq\pcal'\subseteq\pcal$ implies l-isomorphism between
 $\pcal$ and $\pcal'$.

A pre-model $\pcal= (W,\R,\{\lgc w\}_{w\in W},\V)$ is called \textit{bi-finite} if  $W$ and 
$\axioms{\lgc w}$ for every $w\in W$  are finite and $\rules{\lgc w}$ is empty. 
Note that in a bi-finite pre-model, we have  finiteness both as frame 
and logical properties\footnote{A theory  is called finite, if the set of 
its theorems are finite.}.

 We say that a provability model $\pcal= (W,\R,\{\lgc w\}_{w\in W},\V)$ 
 is \textit{decidable}, if   ${\lgc w}\vdash A$ is decidable for every $A$ and $w\in W$. 
 Note that decidability of a finite provability model, implies the decidability of  
 satisfaction relation $\pmodels$.
 
 Remember that a tree, is a relation such that no downward branching happens 
 (see \cref{sec-order}).
\begin{theorem}\label{bases-gl}
For every converse well-founded tree provability pre-model $\pcal$, 
there is a provability model $\bar\pcal$ with following properties:
\begin{itemize}
\item $\bar\pcal$ has necessitation, 
\item $\pcal\subseteq \bar\pcal$\footnote{This means that $\bar{\pcal}$ is also converse well-founded.},
\item For every provability model $\pcal'\supseteq \pcal$ with necessitation, 
we have $\pcal'\supseteq\bar\pcal$. In other words, $\bar\pcal$ is the minimum 
provability model  with necessitation including $\pcal$.
\end{itemize}
Moreover, if $\pcal$ is bi-finite, then $\bar\pcal$ is decidable.
\end{theorem}
\begin{proof}
Let $\pcal=(W,\R,\{\lgc w\}_{w\in W},\V)$. 
We define $\bar\pcal=(W,\R,\{\blgc w\}_{w\in W},\V)$ in which $\blgc w$ is defined 
by recursion on $w\in W$ ordered by $\Rhat$ (see \cref{sec-order}) as follows.  Notice that 
by \Cref{cwf-tree}, $(W,\Rhat)$ is  converse well-founded and hence the recursion is 
a legitimate one indeed. 
As inference rules of $\blgc w$, we only include modus ponens. However we see later in the 
proof that necessitation is admissible to $\blgc w$'s. Thus we may add necessitation 
explicitly to the set of inference rules of $\blgc w$'s without changing the validity of 
formulas.  

Therefore, we assume that 
$\blgc u$ is defined for every $u\hatR w$, and hence 
$\bar\pcal,w\pmodelsp A$   is already defined for every purely modal $A$. 
%(note that here we need the tree structure to have uniqueness of the immediate predecessor).
Define
$$
\axioms{\blgc w}:=\Gamma\cup\Box\Gamma   
\quad\text{and}\quad 
\Gamma:=\axioms{\lgc w}\cup\{A: A \text{ is purely modal and }\bar\pcal,w\pmodelsp A  \}\cup \taut
.
$$ 
In the above formulation, $\Box \Gamma:=\{\Box A: A\in \Gamma\}$ and 
$\taut$ is the  set of all  classical tautologies.  
It is obvious from the definition that $\bar\pcal$ is indeed a provability model.
\\[2mm]
\textit{Claim:} $\blgc w$ includes $\kfour$ and is closed under necessitation.
\\[2mm]
We prove   by induction on $w\in W$ ordered by $\hatR$. As induction hypothesis, 
for every $u\hatR w$  
we have that $\blgc u$ has necessitation  and also includes $\kfour$.
We first show that $\blgc w$ includes $\kfour$.
For this, we need to show that all axioms of $\kfour$
belongs to $\axioms{\blgc w}$. 
By definition of $\blgc w$, all classical tautologies and their boxes are already included 
in $\axioms{\blgc w}$. 
For other axioms $\varphi$ of $\kfour$, it is enough to show 
$\bar\pcal, w\pmodelsp  \varphi$. We have  following cases:
%Hence, it is enough to show that 
%$\bar\pcal, w\pmodelsp  \varphi$ for every axiom $\varphi$ of $\kfour$. 
%However, by induction hypothesis,
%we already have $\bar\pcal, u\pmodels  \varphi$ for every $u\hatR w$. 
%Hence, it is enough to show 
%$\bar\pcal, w'\pmodels  \varphi$ for $w'\Rsim w$ and any  axiom $\varphi$ of $\kfour$:
\begin{itemize}
%\item Classical tautologies: By definition.
%\item $\Box A$ for a classical tautology $A$:   By definition of $\axioms{\blgc w}$.
\item $\Box(A\to B)\to (\Box A\to\Box B)$:  
      By this fact that  $\blgc u$ has modus ponens for every $u\in W$.
\item $\Box A\to \Box\Box A$: Since by induction hypothesis, $\blgc u$ is closed 
under necessitation, for every $u\hatR w$.
\item $\Box A$, with $A$ being any of above two axioms: 
By induction hypothesis for any $u\hatR w$ we have 
$\blgc u\vdash A$.   Hence, by definition we have $\pcal,w\pmodelsp \Box A$.
\end{itemize}
It remains only to show that $\blgc w$ is closed under necessitation. So assume that 
${\blgc w}\vdash A$. Hence there is a finite set $\Delta\subset \Gamma$ such that 
$\Delta,\Box\Delta\vdash A$. This implies that $\kfour\vdash \bigwedge\Boxdot\Delta\to A$.
Thus $\kfour\vdash \bigwedge \Box \Delta\to \Box A$. 
Now, given that $\blgc w$ includes 
$\kfour$ and $\Box \Delta$, 
we have $\blgc w\vdash \Box A$, as desired.
\\[2mm]
Hence  $\bar\pcal$, is indeed a provability model with (admissibility of) necessitation. 
Also, by definition we obviously have $\pcal\subseteq\bar\pcal$.
The minimality of $\bar\pcal$, is immediate by its definition. 

Finally we show that if $\pcal$ is bi-finite, then $\bar{\pcal}$ is decidable. 
We do this by induction on $w\in W$, ordered by $\hatR$, 
that $\blgc w$ is decidable. 
So as induction hypothesis, we assume that for 
every $u\hatR w$ we have the decidability of $\blgc u$. This implies that for every 
purely modal formula $B$, we have the decidability of $\bar\pcal,w\pmodelsp B$. 
Since $\lgc w$ is finite, we may define  $\varphi_w:=\Boxdot\bigwedge\axioms{\lgc w}$. 
For a given $A$, we must decide ${\blgc w} \vdash A$. 
The decision algorithm works as follows:
\begin{itemize}
\item First calculate the purely modal uniform pre-interpolant of the formula 
$\varphi_w\to A$, namely $(\varphi_w\to A)^*$. 
\item Given that $(\varphi_w\to A)^*$ is purely modal, by induction hypothesis, 
decide $\bar\pcal, w\pmodelsp (\varphi_w\to A)^*$ and return the outcome.
\end{itemize}
We will prove that this algorithm is correct. First assume that the algorithm returns yes. 
This means that $\bar\pcal,w\pmodelsp (\varphi_w\to A)^*$. Since $(\varphi_w\to A)^*$
is purely modal, this implies that $(\varphi_w\to A)^*\in\axioms{\blgc w}$. 
Since $\vdash (\varphi_w\to A)^*\to (\varphi_w\to A)$, we get 
${\blgc w}\vdash \varphi_w\to A$. Given that ${\blgc w}\vdash \varphi_w$, 
this implies $ {\blgc w} \vdash A$.

For the other way around, let $  {\blgc w}\vdash A$. This means that there is some 
purely modal formula $B$ such that $\bar\pcal,w\pmodelsp B$ and 
$\Boxdot B\vdash \varphi_w\to A$. This, by the third property of purely modal uniform pre-interpolation, implies 
$\vdash \Boxdot B\to (\varphi_w\to A)^*$. Hence by $\bar\pcal,w\pmodelsp \Boxdot B$ we get 
$\bar\pcal,w\pmodelsp (\varphi_w\to A)^*$. Thus the algorithm returns yes. 
\end{proof}

\begin{definition}\label{finitary}
Above theorem, states that any conversely well-founded tree pre-model can be 
extended to a unique (up to l-isomorphism) minimum provability model with  the same 
frame and with necessitation. 
Such unique provability model is called to be \textit{generated} by the 
pre-model.  We call a provability model to be \textit{finitary}, 
if it is generated by a bi-finite conversely well-founded tree
pre-model.
\end{definition}
Finitary provability models are well-behaved, in the sense that they can 
be \textit{constructed} and their validity are decidable. Mentioned properties, 
making them convenient tools, like usual Kripke models.
In the following theorem, we improve our earlier \Cref{completeness-gl}, 
showing that the class of finitary provability models are complete for $\GL$.

\begin{theorem}\label{completeness-gl-finitary}
$\GL$ is complete for finitary provability models.
\end{theorem}
\begin{proof}
    Let $\GL \nvdash A$. 
	Then by \Cref{finfalsGL}, there is some $n$ such that $\GL\nvdash \Box^n\bot\to A$.
    By completeness of $\GL$ for finite irreflexive transitive tree Kripke models, 
    there is some $\kcal = (W, \R, \V)$ such that $\kcal,w_0\models \Box^n\bot$ and 
    $\kcal, w_0 \nmodels A$ 
    for some $w_0 \in W$.
    Let $X$ be a finite set of representatives of formulas modulo equivalence in 
    $\GL_n$, as promised by \Cref{gl-local-tabularity}. Then let
    \[
    \axioms{\lgc w} := \{B \in X \mid \kcal, w \models \Boxdot B\}
    \quad \text{and}\quad 
    \rules{\lgc w}:=\emptyset
    \quad \text{and}\quad
    \pcal := (W, \R, \{\lgc w\}_{w \in W}, \V).
    \]
    Then let 
    $\bar{\pcal}:=(W, \R, \{\blgc w\}_{w \in W}, \V)$ be the provability model 
    generated by $\pcal$ (hence having necessitation and modus ponens).
    Also, define
    \[
    \axioms{{\lgcp w}} := \{B \mid \kcal, w \models \Boxdot B\}
     \quad \text{and}\quad 
     \rules{\lgcp w}:=\{{\sf mp, nec}\}
     \quad \text{and}\quad 
     \pcal' = (W,\R, \{\lgc w '\}_{w \in W},\V).
    \]
    Notice that then $\pcal'$ is equal to the provability model defined in the proof 
    of Theorem~\ref{lem1} and hence $\pcal'$ is a finite tree conversely-wellfounded 
    (as it is finite, irreflexive and transitive) provability model with necessitation 
    such that for any $v \in W$ and formula $B$
    \[
        \kcal, v \models B \quad \text{iff}\quad \pcal',v\pmodels B.
    \]
     By definition it is obvious that $\pcal \subseteq \pcal'$, 
    so by minimality of $\bar{\pcal}$ we obtain $\bar{\pcal} \subseteq \pcal'$.
    If we want to show that $\bar{\pcal}, w \not\pmodels A$ it suffices to show that 
    $\pcal' \subseteq \bar{\pcal}$, as this would imply the l-isomorphism of 
    $\pcal'$ and $\bar{\pcal}$.
    
    First, let us prove that $\blgc v$ includes $\GL_n$, for any $v \in W^\R$.
    By Theorem~\ref{soundness-gl} we know that $\GL$ is sound for conversely 
    wellfounded provability model having necessitation.
    By Theorem~\ref{bases-gl} we know that $\bar{\pcal}$ is a conversely 
    wellfounded provability model having necessitation, so $\bar{\pcal} \pmodels \GL$.
    This means that for any $\GL\vdash B$ we also have $\GL\vdash \Box B$ and hence 
    $\bar\pcal,w\pmodels \Box B$ for any theorem $B$ of $\GL$ and any $w\in W$. 
    Thus, by definition of validity in $\bar\pcal$, 
    for any $\R$-accessible node $v$ we have ${\blgc v}\vdash B$, as desired.
	Given that $\Box^n\bot$ is purely modal, to show that ${\blgc v}\vdash \Box^n\bot$, 
	it is enough to show that $\bar{\pcal} \pmodels \Box^n \bot$. 
	So let $v \in W$ and pick 
    $u \sqsupset v$, we have to show that $\blgc u \vdash \Box^{n-1} \bot$. 
    For this it suffices to show that $\Box^{n-1} \bot \in \axioms{\lgc u}$, 
    or equivalently, that $\kcal, u \models \Boxdot \Box^{n-1} \bot$.
    Since $\kcal \models \Box^n \bot$ we have that $\kcal, u \models \Box^n \bot$, 
    and since $v \R u$ and $\kcal, v \models \Box^n \bot$ we obtain that 
    $\kcal, u \models \Box^{n-1} \bot$ (where we used that $n > 0$).
%    Then $\bar{\pcal} \models \GL_n$, and then (by modal completeness and closure 
%    under modus ponens) for any $v \in W^\R$ we have that $\GL_n \subseteq \blgc v$.

    Finally, pick an arbitrary $v \in W$  
    aiming to show $\lgcp v \subseteq \blgc v$.
    Assume that $ {\lgcp v}\vdash B$. Since $\axioms{\lgcp v}$ is obviously already closed 
    under necessitation and modus ponens, $ {\lgcp v}\vdash B$ is equivalent to 
    $B\in \axioms{\lgcp v}$. Thus from ${\lgcp v}\vdash B$ we 
    get $\kcal, v \models \Boxdot B$.
    On the other hand, 
    there is a $B' \in X$ such that $\GL_n \vdash B \leftrightarrow B'$, 
    and then also $\GL_n \vdash \Box B \leftrightarrow \Box B'$, which implies 
    $\GL_n \vdash \Boxdot B \leftrightarrow \Boxdot B'$.
    Since $\kcal \models \GL_n$, 
    we have $\kcal, v \models \Boxdot B'$.
    Using that $B' \in X$, by definition of 
    $\lgc v$ this gives $ {\lgc v}\vdash B'$ 
    and then $ {\blgc v}\vdash B'$.
    We previously showed that $\GL_n \subseteq \blgc v$, so 
    $ {\blgc v}\vdash B' \leftrightarrow B$ and we can conclude $ {\blgc v}\vdash B$, 
    as desired.
\end{proof}

\begin{remark}
Provability models, are very much aligned with the 
arithmetical provability interpretations. Therefore, in the case of $\GL$,
it is very natural to wonder how to prove arithmetical 
completeness theorem by means of  provability models. 
The first author of this paper, actually used \cite{PLHA,mojtahedi2021hard} 
provability models to \textit{reduce}
provability logic of Heyting Arithmetic ($\HA$) to its $\Sigma$-provability logic 
 \cite{Sigma.Prov.HA}. 
One may use the similar argument and reduce the provability logic of $\PA$
to the $\Sigma$-provability logic of $\PA$, by means of provability models. 
Nevertheless, the reduction in the case of classical provability logic, unlike the 
intuitionistic case, is already simple and could be carried out easily 
by standard Kripke models \cite{reduction}. 
\end{remark}
\section{Provability models for the Interpretability Logic}\label{sec3}
In this section we study the provability models of  the interpretability logic $\ILM$.
In the interpretability logic, there is a binary modal operator $\rhdi$, with the 
following intended 
arithmetical interpretation  for $A\rhdi B$:
\begin{center}
$\PA+A$ interprets $\PA+B$.
\end{center}
For more details see \cite{Visser-Interpretability}.
The arithmetical completeness of the logic of interpretability was proved in  
in \cite{Shav,Berarducci}. 

Interestingly, interpretability for Peano Arithmetic, is equivalent to 
$\Pi_1$-conservativity. In other words, $\PA+A$ interprets $\PA+B$
iff  every $\Pi_1$-consequence of $\PA+B$ is also $\Pi_1$-consequence of 
$\PA+A$ (see \cite{Berarducci}).
%In terms of $\Sigma_1$-preservativity, it can be rephrased as follows.
%  $\PA+A$ interprets $\PA+B$
%iff for every $S\in \Sigma_1$ from $\PA\vdash S\to \neg B$ we get $\PA\vdash S\to\neg A$.

So %if we take $(A\rhdi B):=(\neg B\rhdi\neg A)$, 
%then 
in terms of propositional language, we have 
$A\rhdi B$ iff for every $\lozenge E$, 
$\sft\vdash B\to\lozenge E$ implies $\sft\vdash A\to\lozenge E$,
for a desired modal theory $\sft$.

This simple observation, gives us the idea that how to define provability models 
for the interpretability logic $\ILM$.

 \subsection*{Axioms and rules of $\ILM$}
%The language $\lcali$ is defined as the boolean language 
%$\{\to,\bot\}$
%together with a binary modal connective $\dhr$. 
%Moreover, we take the set of atomics \textit{finite}.

%Let $\sft$ be a logic%\footnote{A modal theory is a set of propositions $T$ in the language $\lcalb$ or $\lcali$ which is closed under modus ponens and necessitation rules.}
%. 
Given a logic $\sfl$, the interpretability logic over $\sfl$,
$\Pre(\sfl)$, is a modal logic with  modus ponens as its only inference rule and 
following axioms:
\begin{itemize}
%\item  $\Box(A\to B)\to(\Box A\to\Box B)$.
%\item $\Box(\Box A\to A)\to \Box A$.
\item  $ \Box(A\to B)\to(A\rhdi B) $.
\item $ (A\rhdi B\wedge B\rhdi C)\to A\rhdi C $.
\item  $ ( B\rhdi   A\wedge   C\rhdi  A)\to ( B\vee  C)\rhdi   A $.
\item $  \lozenge   A\rhdi  A $.
\item  $ A\rhdi B \to 
\big((\Box C\wedge A)\rhdi(\Box C\wedge B)\big) $. (Montagna's axiom)
\item $\Box A$ for $A$ being any of the above mentioned axioms.
%\item  $A$ / $\Box A$ (necessitation).
\item  All theorems of $\sfl$. 
\footnote{Note that  $\sfl$ might be a modal logic in the language $\lcalb$, 
however here we consider all theorems of $\sfl$ in the bigger language $\lcali$.}
%\item $A$ and $A\to B$ / $B$ (modus ponens).
\end{itemize}
%\end{remark}
Then we define $\ILM:=\Pre(\GL)$. Note that this is equivalent to the 
original axiomatisation of the interpretability logic of Peano Arithmetic (as in \cite{Shav,Berarducci}). 
Probably, the most notable thing missing in the current axiomatization of $\ILM$
is the necessitation rule. However, the necessitation is admissible, 
since we already included the box  of every instance of the axiom schemes and 
with the aid 
of the transitivity axiom schema in $\GL$. 

\begin{remark}
Comparing this system   to the  $\iph$
\cite{IemhoffT,Iemhoff2005}, a candidate for $\Sigma_1$-preservativity logic of Heyting Arithmetic $\HA$, one may observe that two axioms are missing here: 
Disjunction and Visser axioms\footnote{One must be careful about the correspondence between the interpretability language and the language of $\iph$. More precisely, $A\rhdi B$ in one language should be translated as $\neg B\rhdi \neg A$ in the other. Thus, the present disjunction axiom, is actually a conjunction axiom there in $\iph$.}.
\end{remark}

\subsection{Veltman models for $\ILM$}
A Veltman model \cite{dejongh-veltman} is a tuple $\vcal:=(W,\{\pce_w\}_{w\in W},\R,\V)$
with following properties:
\begin{itemize}
    \item $W$ is a nonempty set of worlds.
    \item $\R$ is a binary relation on $W$.
    \item $\pce_w$ is a transitive reflexive relation on the set of $\R$-accessible nodes from $w$, i.e.~$(w{\R}):=\{u\in W: w\R u\}$.
    \item $w\R u\R v$ implies $u\pce_w v$.
    \item $u\pce_w v\R z$ implies $u\R z$.
    \item $ \V $ is a binary relation between worlds and propositional variables.
\end{itemize}
Then the relation $\vcal,w\models A$ is defined inductively. It commutes with boolean connectives and:
$$\vcal,u\models A\rhdi B \quad\text{iff}\quad  
\forall\,v\sqsupset u\ 
\big(
\vcal,v\models A 
\Rightarrow    
(\exists\, w\sce_u v\ \vcal,w\models B)
\big).$$
\begin{remark}
By other properties of Veltman models,
the following definition for validity of $A\rhdi B$ is equivalent to the original one defined above:
$$\vcal,u\models A\dhr B \quad\text{iff}\quad  \forall\,v\sqsupset u\ \big(
(\exists\, w\sce_u v\ \vcal,w\models A)
\Rightarrow 
(\exists\, w\sce_u v\ \vcal,w\models B)
\big).$$
This alternative definition is symmetric and better 
aligned with the provability semantics which we introduce later for $\ILM$.
\end{remark}

\begin{theorem}\label{ILM-models}
    $\ILM$ is sound and complete for Veltman models. Furthermore, 
    the completeness also holds for finite Veltman models.    
\end{theorem}
\begin{proof}
    The soundness is straightforward and left to the reader. For the 
    completeness  see  Appendix A in \cite{Berarducci},  or 
    \cite{Velt-Jongh}. 
\end{proof}
Thanks to the completeness for finite models we have the following corollary.
\begin{corollary}\label{finfals}
$\ILM\vdash A$ iff $\ILM\vdash \Box^k\bot\to A$ for all 
$k\in\nat$.
\end{corollary}

Let us define $\ilmn n:=\ILM+\Box^n\bot$. 
\begin{lemma}\label{ilm-local-tabularity}
For every $n$, $\ilmn n$ is locally tabular, 
i.e.~modulo $\ilmn n$-provable equivalence, there are finitely many 
formulas over the finite language. 
Moreover, One may calculate a finite set $X_n$ includes an $\ilmn n$-equivalent of 
every formula.
\end{lemma}
\begin{proof}
Given that Veltman models of $\ilmn n$ are the models with height $n$, 
one may easily observe that every formula is equivalent to a formula in 
which the maximum number of nesting 
$\rhdi$'s are $n$. Then one may easily construct the desired finite set by induction on $n$.
%Before we continue, let us fix some notation. 
%For a given finite set of formulas $X$, let $\Bool X$ be a finite set 
%such that for every boolean combination  of formulas in $X$, there is a member of 
%$\Bool X$ which is equivalent in classical logic to that formula. 
%It is obvious that for every finite $X$, the finite set $\Bool X$ can be 
%effectively calculated.
%
%We define first the sets $X_n$   by recursion on $n$. 
%\begin{itemize}
%\item $X_0:=\Bool{\atom}$.
%\item $X_{n+1}:=\Bool{X_n\cup\{ A\rhdi B: A,B\in X_n\}}$.
%\end{itemize}
%Then by induction on $n$, we show that every formula $A$  belongs to $X_n$, modulo 
%$\ilmn n $-provable equivalence relation.
%\begin{itemize}
%\item $n=0$. Notice that in $\ilmn 0$, every formula $A\rhdi B$ is equivalent to $\top$. This means that in $\ilmn 0$, every formula, has an $\ilmn 0$-equivalent non-modal formula, i.e.~a boolean combination of atomic formulas, which by definition belongs to 
%$X_0$.
%\item 
%\end{itemize}
\end{proof}

\subsection*{Conservativity}
Given a set  $X$  of formulas and a theory $\sft$, 
define $A\pres X \sft B$ as follows:
$$A\pres X {\sft} B \quad \text{iff}\quad \forall \, E\in X\  (\ \sft\vdash B\to E \Rightarrow 
\ \sft\vdash A\to E).$$
%In the above notation, whenever $Y$ is empty or $X=\{\top\}$,  
%we may omit them in $A\pres X Y B$.
For our future notations, we also define 
$$\Pi:=\{\lozenge E: E\in\lcali\}.$$
In a sense, $\Pi$ is the set of formula in the propositional language which corresponds 
to the set of forst-order $\Pi_1$-formulas.

\subsection{Provability models for $\ILM$}
For a given provability model (pre-model)
$\pcal=(W,\R,\{\lgc w\}_{w\in W},\V)$, we may extend
the validity of formulas to the language with 
a binary modal operator $\rhdi$ as follows:
$$\pcal,w\pmodels A\rhdi B\quad\text{iff}\quad 
\forall u\sqsupset w \ (
 A\pres \Pi {\lgc u} B
).
$$
In other words:
$$\pcal,w\pmodels A\rhdi B\quad\text{iff}\quad 
\forall u\sqsupset w\forall E\ (
\ \vdash_u B \to \lozenge  E
\Rightarrow 
\ \vdash_u A\to \lozenge E
).
$$
Notice that, having this definition for the validity of $A\rhdi B$ and defining
$\Box A:=\neg A\rhdi \bot$, then we have our earlier definition for the validity of 
$\Box A$ for every provability model with necessitation, namely 
$$\pcal,w\pmodels \Box A\quad \text{iff}\quad \forall u\sqsupset w\ (\ \vdash_u A).$$

Also notice that the notion of \textit{purely modal} in the extended language, 
would be extended in the following sense: a formula is purely modal 
(in the language $\lcali$) if it is boolean combination of 
formulas of the form $B\rhdi C$.  

\begin{theorem}\label{Soundness-ILM}
$\ILM$ is sound for conversely well-founded provability models with necessitation.
\end{theorem}
\begin{proof}
Let $\pcal=(W,\R,\{\lgc w\}_{w\in W},\V)$ be a provability model with necessitation. 
Notice that if $A$ is a theorem of $\GL$, then by \Cref{soundness-gl} we have that $\pcal \pmodels A$.
In particular, for any  purely modal axiom $A$ of $\GL$,
we have that $\pcal, w \pmodels^+ A$, 
and since $A$ is purely modal, by modal completeness we get $\vdash_w A$.  
For other axioms  $A$ of $\GL$ which are not purely modal, 
namely the classical tautologies, we have $\vdash_w A$ since $\pcal$ is assumed to be 
classical (as a pre assumption included in the definition of provability model).
This also implies that $\GL \subseteq \lgc w$.

By induction on the complexity of proof $\ILM\vdash \varphi$, 
we show that $\pcal,w\pmodels A$.
We only show the soundness of the axioms.  Soundness of the  modus ponens is trivial.
\begin{itemize}
    \item $A = \Box(B \to C) \to B \rhdi C$: assume $\pcal,w \pmodels \Box(B \to C)$ 
    and hence $\forall u \sqsupset w (\ \vdash_u B \to C)$. 
    Take some  $u \sqsupset w$  and $E$ such that $\vdash_u C\to\lozenge E$. 
    Then by $\vdash_u B \to C$ we may infer 
    (using classical logic which is included in $\lgc u$) 
    that $\vdash_u B\to\lozenge E$, as desired.
    \item $A = (B \rhdi C \wedge C \rhdi D) \to B \rhdi D$: 
    assume $\pcal, w \pmodels B \rhdi C$ and $\pcal, w \pmodels C \rhdi D$, 
    pick $u \sqsupset w$ and $E$ such that $ \vdash_u D\to \lozenge E$. We need 
    to show $\vdash_u B\to\lozenge E$.
    Since $\pcal, w \pmodels C \rhdi D$, $ \vdash_u D\to \lozenge E$ 
    implies $ \vdash_u C\to \lozenge E$. 
    Since $\pcal, w \pmodels B \rhdi C$, we obtain $\vdash_u B\to\lozenge E$, as desired.
    \item $A = (C \rhdi B \wedge D \rhdi B) \to (C\vee D) \rhdi B$: 
    assume that $\pcal, w \pmodels C \rhdi B$ and $\pcal, w \pmodels D \rhdi B$. 
    Pick arbitrary $u \sqsupset w$ and $E$ such that $ \vdash_u B\to\lozenge E$. 
    We need to show  $ \vdash_u  (C \vee D)\to \lozenge E$.
    By $\pcal, w \pmodels C \rhdi B$, we have  $ \vdash_u  C\to\lozenge E$. 
    Also  by $\pcal, w \pmodels D \rhdi B$ we have $  \vdash_u  D\to\lozenge E$.
    By classical propositional reasoning we obtain $ \vdash_u  (C \vee D)\to \lozenge E$,
    as desired.
    \item $A = \lozenge B \rhdi B$: let $u \sqsupset w$ and $E$ 
    such that $  \vdash_u B\to\lozenge E$.
    Since $\lgc u$ is closed under necessitation, modus ponens and have the axioms of 
    $\GL$ we have that $  \vdash_u  \lozenge B\to\lozenge\lozenge E$.
    Also $ \vdash_u \lozenge\lozenge E \to \lozenge E$, 
    so we can conclude the desired $  \vdash_u \lozenge B\to \lozenge E$.
    \item $A = (B \rhdi C) \to ((\Box D \wedge B) \rhdi (\Box D \wedge C))$: 
    assume that $\pcal, w \pmodels B \rhdi C$ and pick $u \sqsupset w$ and $E$ 
    such that $  \vdash_u (\Box D \wedge C)\to\lozenge E$.
    We need to show that $\vdash_u (\Box D\wedge B)\to\lozenge E$.
    By classical propositional reasoning, we have that 
    $ \vdash_u C\to (\lozenge E\vee\lozenge \neg D)$. 
    Since $\lgc u$ contains $\GL$ we have that 
    $\vdash_u (\lozenge E \vee \lozenge \neg D) \Leftrightarrow \lozenge (E \vee \neg D)$. 
    Hence $\vdash_u C\to \lozenge (E \vee\neg  D) $.
    Using $\pcal,w \pmodels B \rhdi C$ we obtain that 
    $ \vdash_u B\to \lozenge  (E \vee\neg D) $. 
    Thus  we can conclude  that 
    $ \vdash_u (B\wedge\Box D)\to\lozenge E$, as desired.
    \item $A = \Box B$ for $B$ being one of the other axioms of $\ILM$: 
    we need
the induction hypothesis on $u\sqsupset w$ and the modal completeness property.
    \item $A$ is a theorem of $\GL$: argued at the beginning of this proof.
    \qedhere
\end{itemize}
\end{proof}

To enable ourselves to prove completeness of $\ILM$ for provability models, 
we need to refine the structure of Veltman models by unravelling. 
Intuitively, we need to get rid of the the dependency of $\pce_w$ on $w$, for which we use the famous method of unravelling:
\begin{definition}\label{def-unrav-Veltman}
Let $\vcal:=(W,\{\pce_w\}_{w\in W},\R,\V)$ be a finite Veltman model. We define 
 $\vcal_0:=(W_0,\R_0,\pce_0,\V_0)$ as follows:
 (see \cref{Fin-seq} for the notations on finite sequences.)
 \begin{itemize}
\item $W_0:=$ The set of all 
finite non-empty sequences $\sigma:=\langle w_1,\ldots,w_n\rangle$ such that $w_1\R w_2\ldots \R w_n$. 
\item $\sigma\R_0 \tau$ iff $\sigma=\iseg\tau$.
In other words, $\sigma\R_0 \tau$ iff 
there is some $w\in W$ s.t.~$\tau=\sigma*\langle w\rangle$.
\item $\sigma\pce_0 \tau$ iff there is some $\eta\in W_0$ such that 
we have $\eta\R_0\sigma,\tau$ and 
$\fel{\sigma}\pce_{\fel{\eta}} \fel{\tau}$.
In other words, $\sigma\pce_0 \tau$ iff
there are some $w,u,v\in W$ and $\eta\in W_0	$ such that 
$\sigma=\eta*\langle u\rangle$ and $\tau=\eta*\langle v\rangle$ and $w=\fel \eta$ and $u\pce_w v$.
\item $\sigma\V_0 p$ iff $\fel\sigma\V p$.
\end{itemize}
Observe that $\vcal_0$ as defined  above, is  a tree.
Although $\vcal_0$ is not a Veltman model, still we may define
validity of formulas $B\in\lcali$ for it, in the same way as we did for $\ILM$. 
More precisely, we define 
$$\vcal_0,\sigma\models C\rhd D \quad\text{iff}\quad  
\forall\,\tau\sqsupset_0 \sigma\ 
\big((\exists\, \eta\sce_0 \tau\ \vcal_0,\eta\models C)\Rightarrow 
(\exists\, \eta\sce_0 \tau\ \vcal_0,\eta\models D)\big).$$ 
\end{definition}

Then, as expected, we have the correspondence of the  Veletman model with its unravelling:

\begin{lemma}\label{lem-unrav}
For every $\sigma\in W_0$ and $B\in\lcali$ we have 
$$\vcal_0,\sigma\models B\quad \text{iff}\quad \vcal,\fel\sigma\models B.$$
\end{lemma}
\begin{proof}
Easy induction on the complexity of 
$B\in\lcali$.
\end{proof}

\begin{remark}\label{lem-unrav2}
By above lemma \ref{lem-unrav} we have 
$$\vcal_0,\sigma\models C\rhdi D \quad\text{iff}\quad  
\forall\,\tau\sqsupset_0 \sigma\ 
\big(\vcal_0,\eta\models C\Rightarrow 
 (\exists\, \eta\sce_0 \tau\ \vcal_0,\tau\models D) \big).$$ 
\end{remark}

\begin{theorem}\label{completeness-ILM}
$\ILM$ is complete for finite tree conversely well-founded provability models with necessitation.
\end{theorem}
\begin{proof}
Let $\ILM\nvdash A$. 
Then by completeness if $\ILM$ for Veltman models
(\Cref{ILM-models}), there is a finite Veltman model 
$\vcal:=(W,\{\pce_w\}_{w\in W},\R,\V)$ such that 
$\vcal,w\nmodels A$ for some $w\in W$.
Let 
$\vcal_0:=(W_0,\R_0,\pce_0,\V_0)$ as defined in \Cref{def-unrav-Veltman}.

Then define 
$\pcal:=(W_0,\R_0,\{\lgc \sigma\}_{\sigma\in W_0},\V_0)$ as follows. 
We only need to define $\lgc \sigma$ for $\R_0$-accessible 
nodes $\sigma$:
$$\axioms{\lgc \sigma}:=\{B\in\lcali: \forall\, \tau\sce_0 \sigma\ (\vcal_0,\tau\models B)\}
\quad \text{and} \quad 
\rules{\lgc\sigma}:=\emptyset
.$$
\textit{Claim 1.}   $\pcal,\sigma\pmodels B$ iff $\vcal,\fel\sigma\models B$ iff 
$\vcal_0,\sigma\models B$.
\\[2mm]
\textit{Proof of the claim 1.}
First notice that by \Cref{lem-unrav}, we already have equivalence of  the two latter
statements. For the equivalence of first two statements, we use induction on $B\in\lcali$. 
All cases are obvious, except $B=C\rhdi D$, for which we have following argument. 
First assume that $\pcal,\sigma\pmodels C\rhd D$ seeking to show 
$\vcal_0, \sigma\models C\rhd D$. For this end, let $\tau\sqsupset_0 \sigma$ 
such that  $\vcal_0,\tau\models C$. We need to show 
$\vcal_0,\eta\models D$ for some $\eta\sce_0 \tau$.
By definition of 
$\lgc \tau$, our assumption implies $\neg C\nin\axioms{\lgc\tau}$. 
Thus by $\pcal,\sigma\pmodels  C\rhd D$, we get 
$\neg D\nin\axioms{\lgc\tau}$\footnote{Notice that here we are using this fact that 
$\lgc \tau$ does not have any inference rules, and hence ${\lgc\tau}\vdash\varphi$ 
is equivalent to $\varphi\in\axioms{\lgc\tau}$}. 
Again, by definition of $\axioms{\lgc \tau}$ we get 
$\vcal_0,\eta\models D$ for some $\eta\sce_0\tau$, as desired. 

Next assume that $\vcal_0, \sigma\models C\rhd D$ seeking to show 
$\pcal,\sigma\pmodels C\rhd D$. Let $\tau\sqsupset_0 \sigma$ such that 
$\nvdash_{\tau}C\to \lozenge E$. 
We must show $\nvdash_{\tau} D\to\lozenge E$. 
Since $\vcal_0$ and $\vcal$ are equivalent, by validity of the Montagna's axiom in $\vcal$,
we get its validity in $\vcal_0$ as well. This means that by our assumption, we also 
have  $\vcal_0, \sigma\models (\Box \neg E\wedge C)\rhd(\Box \neg E\wedge D)$.
By $C\to\lozenge E\nin\axioms{\lgc\tau}$, for some $\eta\sce_0\tau$ we have 
$\vcal_0,\eta\models C\wedge\Box \neg E$.  Therefore,  
$\vcal_0, \sigma\models (\Box \neg E\wedge C)\rhd(\Box \neg E\wedge D)$ implies   
$\vcal_0,\eta\models \Box \neg E\wedge D$ for some $\eta\sce_0\tau$.
Thus $(D\to\lozenge E)\nin \axioms{\lgc \tau}$, as desired. 
This finishes the proof of claim 1.
\\[2mm]
\textit{Claim 2.}   $\lgc\sigma$  includes classical validities, and modus ponens and necessitation are admissible to it. Moreover it satisfies modal completeness.
\\[2mm]
\textit{Proof of the claim 2.}
The fact that $\lgc\sigma$ includes all classical tautologies and is closed under 
modus ponens, are immediate by definition of $\lgc \sigma$. Also 
the closure under necessitation,   
is a direct consequence of claim 1 and the soundness of the axiom 
 $\lozenge B\rhdi B$ in Veltman models.
 
It only remains to show that it satisfies 
modal completeness.  Let $\pcal,\sigma\pmodelsp B$ for some purely modal $B$ and some 
$\R_0$-accessible $\sigma$. We need to show $B\in\axioms{\lgc\sigma}$.
Since $\sigma$ is $\R_0$-accessible, there is a unique $\tau\in W_0$ and a unique 
$w\in W$ s.t.~$\sigma=\tau*\langle w\rangle$. 
Then by $\pcal,\sigma\pmodelsp B$, for every $\eta\sqsupset_0\tau$ we have 
$\pcal,\eta\pmodels B$. On the other hand, $\sigma\sce_0 \sigma'$ implies  
$\tau\sqsubset_0 \sigma'$. Therefore, for every $\sigma'\sce_0\sigma$ we have 
$\pcal,\sigma'\pmodels B$. Thus $B\in\axioms{\lgc\sigma}$, as 
desired\footnote{Notice that this even holds for every $B$, 
and not only $B$'s which are purely modal.}.  
\\[2mm]
Now with the claim 1, we have $\pcal\npmodels A$. With the claim 2, we have that $\pcal$
is a finite, conversely well-founded provability model with necessitation, as desired.
\end{proof}

\begin{remark}\label{remark-1}
The motivation to define 
$w\pmodelsp A$ such that it consider validity of $A$ in the model generated by 
a predecessor of $w$ (excluding the predecessor itself), is quite technical 
for the aim of proving the modal completeness, in the previous theorem. 
Otherwise, we could have taken the following definition for $w\pmodelsp A$ which sounds 
more natural and easier to work with: $w\pmodelsp A$ iff  $u\pmodels A$ for any $u$ 
 accessible from $w\sqsupseteq^*u$, in which $\sqsubseteq^*$ is the transitive and 
 reflexive closure of $\R$.
 Actually, for the case of $\GL$, we could have taken this alternative requirement in the 
 definition of provability models, and the soundness and completeness theorems 
 would have worked. 
\end{remark}

\subsection{Decidability of Provability models for $\ILM$}
In this subsection,  we show that for $\ILM$ also one may construct 
independent new provability models and moreover, in case of finite models, we have 
decidability of validity in the provability model. Notice here that, the decidability 
of validity in provability models in the language $\lcali$, is not merely given to us by 
decidability of the theories and finiteness of the frame. This is because for  
$w\pmodels A\rhdi B$, we need to check $\Pi$-conservativity, 
which asks for verifying the following for \textit{infinitely} many $\lozenge E$'s:
$$\forall u\sqsupset w\forall E\ (
\ \vdash_u B \to \lozenge  E
\Rightarrow 
\ \vdash_u A\to \lozenge E ).$$

The only difference between the following theorem and \Cref{bases-gl} is  on the language. 
In the following theorem, the language is $\lcali$, and hence the set of purely modal
formulas includes more formulas than the one for the language $\lcalb$. As a result, 
the notion of modal completeness for provability models of the two languages are
different. 
%That's why we need to check that everything works here for the language $\lcali$ too. 
%Nevertheless, we skip the ide

\begin{theorem}\label{bases-ilm}
For every converse well-founded tree provability pre-model $\pcal$, 
there is a provability model $\bar\pcal$ with following properties:
\begin{itemize}
\item $\bar\pcal$ has necessitation, 
\item $\pcal\subseteq \bar\pcal$\footnote{This means that $\bar{\pcal}$ is also converse well-founded.},
\item For every provability model $\pcal'\supseteq \pcal$ with necessitation, 
we have $\pcal'\supseteq\bar\pcal$. In other words, $\bar\pcal$ is the minimum 
provability model  with necessitation including $\pcal$.
\end{itemize}
Moreover, if $\pcal$ is bi-finite, then $\bar\pcal$ is decidable.
\end{theorem}
\begin{proof}

The same proof of \Cref{bases-gl} works here as well, 
except for the statement regarding its decidability, for which we have following argument.

We show that if $\pcal$ is bi-finite, then $\bar{\pcal}$ is decidable. 
First, we show by induction on $w\in W$ ordered by $\hatR$ that $\blgc w$ is decidable.
The same procedure provided in the proof of \Cref{bases-gl} (adapting the purely modal uniform pre-interpolation to the $\lcali$ language, which is straightforward) works.
However, this does not directly imply the decidability of $\bar{\pcal}, w \pmodels A$ 
as in the previous case.
The problem being that, in principle, to check if $\bar{\pcal},w \pmodels A \rhd B$ 
we need to check for each $u \sqsupset w$ (this is okay since the number of such 
worlds is finite) if for any modal formula $E$ we have that 
${\blgc u}\vdash B \to \lozenge E$ implies $\blgc u \vdash A\to\lozenge E$.
So we need to limit the possible choices of $E$ or otherwise this process is not feasible.

First, note that since $\bar\pcal$ is with finite conversely well-founded frame, 
there is some $n$ such that $\bar\pcal\pmodels \Box^n\bot$, and given the model this 
$n$ can be computed (for example, we can pick $n$ to be the number of worlds).
Then by \Cref{ilm-local-tabularity} we have some finite set $X$ 
of formulas such that every  formula belongs to $X$, modulo 
$\ilmn n$-provable equivalence relation.
Then, to decide if $\bar{\pcal},w \pmodels A \rhd B$ go through each $u \sqsupset w$ 
and each $E \in X$, decide if $\blgc u \nvdash B\to\lozenge E$ or 
$\blgc u \vdash A\to\lozenge E$.
If the answer is affirmative for each $u \sqsupset w$ and each $E \in X$, then return yes,
otherwise return no.
Since $\bar{\pcal} \pmodels \Box^n \bot$ then $\lgc u \vdash \ilmn n$, 
from which the correctness of the algorithm trivially follows.

For the other kinds of formulas decidability is managed as usual, just deciding 
first smaller subformulas.
%So as induction hypothesis, we assume that for 
%every $u\hatR w$ we have the decidability of $\blgc u$. 

% Then the same algorithm and the same correctness of the algorithm which was provided in 
% the proof of \Cref{bases-gl} works here, except for this difficulty to show that 
% Induction hypothesis implies that for every 
%purely modal formula $B$, we have the decidability of $\bar\pcal,w\pmodelsp B$, for which we have the following argument. 

%  First, note that since $\bar\pcal$ is with finite conversely well-founded frame, 
%  there is some $n$ such that $\bar\pcal\pmodels \Box^n\bot$. 
%  Then by \Cref{ilm-local-tabularity} we have some finite set $X$ 
%  of formulas such that every  formula belongs to $X$, modulo 
%  $\ilmn n$-provable equivalence relation. 
  
%  Next we show by induction on $B$ that $\bar\pcal,w\pmodelsp B$, given the assumption that 
%  $A\in\blgc u$ is decidable for every $A$ and $u\hatR w$. All cases are obvious, 
%  except $B=C\rhd D$. Let $w_0\R w$ and $v\sqsupset u\sqsupset w_0$. We need to decide 
%  ``for all $E\in\lcali$, $(\Box E\to C)\in\blgc v$ implies $(\Box E\to D)\in \blgc v$". 
%  Given that $\blgc v\vdash \ilmn n$, we only need to decide the mentioned statement 
%  for finitely many $E\in X$, which can be effectively done since $w\Rhat v$ and hence 
%  by induction hypothesis $\blgc v$ is decidable.
\end{proof}

Recall from \Cref{finitary} that 
  a provability model is called  finitary, 
if it is generated by a bi-finite conversely well-founded tree
pre-model.

\begin{theorem}\label{completeness-ilm-finitary}
$\ILM$ is complete for finitary provability models.
\end{theorem}
\begin{proof}
Let $\ILM\nvdash A$. Then by \Cref{finfals}, there is some $n\in\nat$ such that 
$\ILM\nvdash\Box^n\bot\to A$. 
By completeness of $\ILM$ for finite converse well-founded 
Veltman models (\Cref{ILM-models}), 
there is some $\vcal:=(W,\{\pce_w\}_{w\in W},\R,\V)$ 
such that $\vcal\nmodels A$ and $\vcal\models \Box^n\bot$.  
Then let $\vcal_0:=(W_0,\R_0,\pce_0,\V_0)$ be as in \Cref{def-unrav-Veltman}.
%Let $\sub{A}$ be the set of all sub-formulas of $A$, including $A$ itself.
Also let $X$ be the finite set of all formulas modulo $\ilmn n$-provable equivalence, 
as given by \Cref{ilm-local-tabularity}.  
Define 
$$\axioms{\lgc \sigma}:=\{B\in X: \forall\, \tau\sce_0 \sigma\ (\vcal_0,\tau\models B)\}
\quad\text{and}\quad 
\rules{\lgc \sigma}:=\emptyset
.$$
%$\lgc \sigma:= \{B\in  X: \vcal_0,\sigma\modelsp B\}$.
Finally define the pre-model 
$\pcal:=(W_0,\R_0,\{\lgc \sigma\}_{\sigma\in W^{\R_0}_0},\V_0)$ and let $\bar\pcal$ 
be the provability model generated by $\pcal$. 
Let also $\pcal'$ defined as $\pcal':=(W_0,\R_0,\{\lgc \sigma '\}_{\sigma\in W^{\R_0}_0},\V_0)$ in which 
$$\axioms{\lgcp \sigma}:=\{B\in \lcali: \forall\, \tau\sce_0 \sigma\ (\vcal_0,\tau\models B)\}
\quad\text{and}\quad 
\rules{\lgcp \sigma}:=\emptyset
.$$
 Notice that $\pcal'$ as defined above, is equal to the provability model defined in the proof of \Cref{completeness-ILM} and hence $\pcal'$ is indeed (l-isomorphic to) 
 a finite tree conversely well-founded provability model with necessitation. 
 Moreover, obviously we have $\pcal\subseteq\pcal'$. Hence by minimality of $\bar\pcal$ 
 (\Cref{bases-ilm}), we have $\bar\pcal\subseteq\pcal'$. 
 
Observe that since $\vcal\models \ilmn n$, 
by \Cref{Soundness-ILM,bases-ilm} we have  
$\bar\pcal\pmodels \ilmn n$ and also $\blgc\sigma$ includes $\ilmn n$ for every $\sigma$.  
This implies that $\blgc\sigma\supseteq\lgc \sigma '$ and hence $\pcal'\subseteq\bar\pcal$.
Together with $\bar\pcal\subseteq\pcal'$, this implies that 
$\bar\pcal$ is l-isomorphic to $\pcal'$. 
Now given that by proof of \Cref{completeness-ILM} we have $\pcal'\npmodels A$, we are done.
\end{proof}

\begin{question}\label{q-completeness}
Is it possible to simplify or having an alternative proof of arithmetical completeness
of $\ILM$, by means of provability models, instead of Veltman models?
\end{question}
\section{Provability models for the Poly-modal Provability Logic}\label{sec-GLP}
The poly-modal provability logic $\GLP$, is a modal logic in the language $\lcalo$. 
The introduction of $\GLP$, goes back to \cite{japaridze1985polymodal}, 
in which he proves the soundness and completeness of $\GLP$ for arithmetical 
interpretations, in which $\boxn n $ is interpreted as provability predicate for 
$\PA$ plus all true $\Pi_n$-sentences.  
After that, L.~Beklemishev brought $\GLP$ in to attentions for other reasons, namely 
by finding some interesting proof-theoretic applications
(see \cite{BEKLEMISHEV2004103}).

Nevertheless, it is known that polymodal provability logic $\GLP$, 
is incomplete for any class of Kripke frames. 
The reason is that the natural frame properties corresponding to the axiom schemes are 
conflicting. However, there are several subsystems of $\GLP$ for which there are
completeness for Kripke frames \cite{ignatiev1993strong,beklemishev2011simplified}. 

On the other hand, we  have topological semantics for $\GLP$ 
\cite{BEKLEMISHEV20131201,Fernández-Duque_Joosten_2013} and also some other infinite 
Kripke-style models \cite{BEKLEMISHEV2010756}.

In this section, we consider provability models for 
$\GLP$. We prove soundness and completeness of $\GLP$ for 
a natural class of provability models. 
However, at leas at this early stage, we are not able to prove decidability 
and constructibility of such models of $\GLP$. More precisely, we only can 
prove completeness via  the canonical model construction.  
We keep this as an open question, that 
if we can find a class of decidable provability 
models\footnote{Provability models in which the validity of formulas at nodes are
decidable.}, for which we have soundness and completeness of $\GLP$.

\begin{definition}\label{poly-p-m}
A poly-provability pre-model (pre-model if no confusion is likely)
is a tuple $\pcal=(W,\{\R_n\}_{n\in\nat}, \{\plgc w n\}_{n\in\nat, w\in W},\V)$ with following properties:
\begin{itemize}
\item $W$ is a non-empty set of worlds.
\item $\plgc w n$ is a theory (in the sense that defined in \cref{logicsec}) 
for every $n\in\nat$ and $w\in W^{\R_0}$.
\item $\R_n$ is a binary relation on $W$ for every $n\in\nat$.
\item ${\V}\subseteq (W\times \atom)$ is the valuation for variables.
\end{itemize}
Then we define $\pcal,w\pmodels A$ by induction on $A$
such that it commutes with boolean connectives and 
$\pcal,w\pmodels p$ iff $w\V p$ and finally:
$$\pcal,w\pmodels \boxn n B\quad \text{iff}
\quad 
\forall\, u\sqsupset_n w\ (\plgc u n \vdash B). 
$$
%Let $\R$ defined as $\R_0$. 
Define  $\pcal,w\pmodelsp_n A$ iff 
for any 
$u\sqsupseteq_n^+ w$ we have $\pcal,u\pmodels A$.
Also we define ${\R}:={\R_0}$.
We say that  $\pcal$ 
is a \textit{poly-provability model} iff for every $n$
and any $w\in W^{\R}$ we have:
\begin{itemize}
\item $\plgc w n $ is classical, i.e.~all classical tautologies are $\plgc w n$-derivable 
and moreover, ${\sf mp}\in\rules{\plgc w n}$.
%\item $\plgc w n$ is closed under modus ponens, 
%i.e.~$A\in\plgc w n$ and $A\to B\in\plgc w n$ 
%implies $B\in\plgc w n$. 
\item Modal completeness: If $A $ is purely modal 
%and $u\R_n w$ 
and $\pcal,w\pmodelsp_0 A$ then $\plgc w 0\vdash A$.
\end{itemize}
Moreover, if $\pcal$ satisfies following properties, 
we call it  poly-provability $\GLP$-model:
\begin{itemize}
\item Poly-necessitation:  
$\pnec n\in\rules{\plgc w n}$, in which  
\begin{prooftree}
\Ax{$A$}\RLa{.}
\LLa{$\pnec n$}
\UI{$\boxn n A$}
\end{prooftree}
\item Poly-L\"ob's rule: 
$\plob n\in\rules{\plgc w n}$, in which  
\begin{prooftree}
\Ax{$\boxn n A\to A$}\RLa{.}
\LLa{$\plob n$}
\UI{$A$}
\end{prooftree}
%%$u\R_n w$ and 
%$(\boxn n A\to A)\in\plgc w n$ implies 
%$A\in\plgc w n$.
\item Ascending: ${\R_{n+1}}\subseteq {\R_n}$ 
and $\plgc w n \subseteq\plgc w{n+1}$.
\item $\Pi$-completeness: 
if $u\R_{n+1} w$ and 
$\pcal,u\pmodels \neg\boxn n A$, then  $\plgc w {n+1}\vdash \neg\boxn n A$.
\end{itemize}
\end{definition}
$\GLP$ is the logic with following axioms and rule:
\begin{itemize}
\item All classical tautologies.
\item $\boxn n (A\to B)\to(\boxn n A\to\boxn n B)$.
\item $\boxn n A\to\boxn n \boxn n A$.
\item $\boxn n (\boxn n A\to A)\to\boxn n A$.
\item $\boxn n A\to\boxn {n+1} A$.
\item $\neg\boxn n A\to\boxn {n+1}\neg\boxn n A$.
\item $\boxn 0 A$ for $A$ being any of the above axioms.
\item Modus ponens.
\end{itemize}
Note that modus ponens is the only inference rule of 
$\GLP$ and we don't have necessitation in our formulation of $\GLP$. 
However, as a set of theorems, this formulation is equal to the original 
axiomatization of $\GLP$. 

%Let us define \textit{poly-provability $\GLP$-model} to be 
%a poly-provability model with following properties:
%\begin{itemize}
%\item Poly-necessitation.
%\item Poly-L\"ob's rule.
%\item $\Pi$-complete.
%\item Ascending.
%\end{itemize}
\begin{theorem}\label{GLP-sound}
$\GLP$ is sound for poly-provability $\GLP$-models.
\end{theorem}
\begin{proof}
Let $\pcal=(W,\{\R_n\}_{n\in\nat}, \{\plgc w n\}_{w\in W},\V)$ be a poly-provability model 
with mentioned properties. We only need to show that 
all axioms of $\GLP$ are valid at $\pcal$. The 
soundness of all axioms
of $\GL$ (for the $\boxn n $ instead of $\Box$)
can be proved in the same way we did in 
\Cref{soundness-gl}. Nevertheless we give a hint also 
here how each axiom can be treated:
\begin{itemize}
\item $\boxn n (A\to B)\to(\boxn n A\to\boxn n B)$: 
By having modus ponens for every
$\plgc w n$.
\item $\boxn n A\to\boxn n \boxn n A$: 
By poly-necessitation.
\item $\boxn n (\boxn n A\to A)\to\boxn n A$: 
By poly-L\"ob's rule.
\item $\boxn n A\to \boxn {n+1} A$: 
Since $\pcal$ is ascending.
\item $\neg\boxn n A\to \boxn {n+1}\neg\boxn n A$:
By $\Pi$-completeness.
\item $\boxn 0 A$ for $A$ being any of above 
axioms: By modal completeness.
\item $A$ is a classical tautology: Obvious.
\item $\boxn 0 A$ for $A$ being a classical tautology:
Since each $\plgc w 0$ includes classical tautologies.\qedhere
\end{itemize}
\end{proof}

\begin{theorem}\label{GLP-completeness}
$\GLP$ is strongly complete for 
 poly-provability $\GLP$-models.
\end{theorem}
\begin{proof}
Let $\Gamma$ be an arbitrary set of formulas and 
$\Gamma\nvdash_\GLP A$. We will define 
a canonical poly-provability model with required properties as follows. 
\begin{itemize}
\item $W$ is the family of all sets $w $ of formulas
which are containing all theorems of 
$\GLP $ and maximally consistent with classical logic.
\item Define $w\R_n u$ iff 
$w,u\in W$ and for any $\boxn n A\in w$ we have 
$A\in u$.
\item Let $\axioms{\plgc w n}:=\{B: B\wedge\boxn n B\in w\}$ 
and $\rules{\plgc w n}:=\emptyset$.
\item $w\V p $ iff $p\in w$.
\item $\pcal:=(W,\{\R_n\}_{n\in\nat},\{\plgc w n\}^{n\in\nat}_{w\in W},\V)$.
%\item Let $W$ be the set of all finite non-empty 
%sequences $\sigma=\langle w_0,k_1,
%w_1, \ldots,w_{n-1},k_n, w_n\rangle$ such that $w_{i-1}\subset_{i} w_i$  for every $0<i\leq n$.
%%and $k_i\leq k_{i+1}$ for every $0<i < n$.
%\item For $\sigma,\tau\in W$, we define 
%$\sigma\R_n \tau$ iff  $\tau=\sigma*\langle k, w\rangle$
%for some  $w\in W_0$ and $k\geq n$.
%\item $\rank \sigma$ for $\sigma\in W$ is defined as 
%the final natural number in the sequence $\sigma$, if there exists any. If there is no such number, we define it as $0$.
%\item $\npred n \sigma$ is defined as the 
%maximum\footnote{This maximality is with respect to the transitive reflexive closure of $\R$, i.e.~$\sqsubseteq^+$.}  
%$\tau\sqsubseteq^+\sigma$ with $\rank\tau\leq n$.
%Note that since $0\leq n $, always there can be found some $\tau\sqsubseteq^+\sigma$ with $\rank\tau\leq n$.
%Thus, a maximal one also exists. Also notice that if we already have $\rank \sigma\leq n$, then 
%$\npred n \sigma=\sigma$.
%\item Let  $n:=\rank \sigma$ and define
%$
%\plgc \sigma n:=
%\{B\in\lcalo: \forall \tau\sqsupseteq^+_n \sigma
%\ (B\in \fel\tau)\}. 
%$
%\item $\sigma\V p$ iff $p\in\fel\sigma$.
%\item  
%Finally define $\pcal:=(W,\R, \{\plgc w n\}_{w\in W^\R}^{n\in\nat},\V)$.
\end{itemize}
\textit{Claim 0.}  If for all  $u\sqsupset_n w$ we have 
$C\in u$, then $\boxn n C\in w$.
\\[2mm]
\textit{Proof of the claim 0.} 
Let $\boxn n C\nin w$. We need to find some $u\sqsupset_n w$ such that $C\nin u$.
Define $u_0:=\{E\in\lcalo: \boxn n E\in w\} $ and take $u\supseteq u_0$
a maximal set such that $u\nvdash C$. Note that since $u_0\nvdash C$, 
such maximal set exists. Furthermore, it is not difficult to observe that such $u$
is actually maximal consistent set and hence $u\in W$.  By definition we also have 
$w\R_n u$.
\\[2mm]
\textit{Claim 1.} $\pcal,w\pmodels B$ iff $B\in w$ for every $B\in\lcalo$ and $w\in W$.
\\[2mm]
\textit{Proof of the claim 1.} We prove this by induction on $B$. 
All cases are obvious except for $B=\boxn n C$. 
First let $\boxn n C\in w$. We want to show 
$\pcal,w\pmodels \boxn n C$.  So let $u\sqsupset_n w$.
Since $\boxn n C\in w$, we also have $\boxn n(C\wedge\boxn n C)\in w$. 
Hence by definition of $w\R_n u$, we get $C\wedge\boxn n C\in u$. 
Hence $\plgc u n\vdash C$. 
Thus we may conclude $\pcal,w\pmodels \boxn n C$, as desired.

For the other way around, let $\boxn n C\nin w$. By claim 0, there is some $u\sqsupset_n w$
such that $C\nin u$.
Given that $\plgc u n\subseteq u$ and $C\nin u$, 
we get $ \plgc u n\nvdash C$. 
Thus we may conclude $\pcal,w\npmodels \boxn n C$, 
as desired.
\\[2mm]
\textit{Claim 2.} $\pcal$ is 
a poly-provability model (modulo l-isomorphism).
\\[2mm]
\textit{Proof of the claim 2.}
Given that any $w$ includes $\GLP$, we have $E\wedge \boxn 0 E\in w$ for any classical 
tautology $E$. This implies that $\axioms{\plgc w n}$ includes classical tautologies. 
Since $w$ is closed under modus ponens and $\boxn n(E\to F)\to(\boxn n E\to\boxn n F)
\in w$, we may infer that $\axioms{\plgc w n}$ is closed under modus ponens. 
It only remains to show the modal completeness, for which the argument follows.
Let $B$ is an arbitrary formula such that $\pcal,w\pmodelsp B$. It is enough to show that 
$B,\boxn 0 B\in w$.  Since $\pcal,w\pmodelsp B$ implies $\pcal,w\pmodels B$, 
by claim 1 we have $B\in w$. On the other hand,  $\pcal,w\pmodelsp B$ implies that 
for all $u\sqsupset w$ we have $\pcal,u\pmodels B$. 
Hence by claim 1 and  
claim 0, we have $ \boxn 0 B\in w$. This finishes showing $B,\boxn 0 B\in w$, as desired.
\\[2mm]
\textit{Claim 3.} $\pcal$ is ascending and $\Pi$-complete and 
has poly-necessitation and 
poly-L\"ob's rule (modulo l-isomorphism).
\\[2mm] 
\textit{Proof of the claim 3.}  
\begin{itemize}
\item Ascending: First notice that since $\boxn n B\to\boxn {n+1} B\in w$, we get 
$\plgc w n\subseteq\plgc w {n+1}$.
On the other hand, let $w\R_{n+1} u$ and $\boxn n B\in w$. 
Then since $\boxn n B\to\boxn {n+1} B\in w$, we get $\boxn {n+1} B\in w$ 
and hence by  $w\R_{n+1} u$ we get $B\in u$. Thus we may conclude $w\R_{n} u$, as desired.
\item $\Pi$-complete: Let $\pcal,u\pmodels\neg\boxn n B$ and $u\R_{n+1} w$. We need to show 
$\plgc w {n+1}\vdash\neg\boxn n B$. By $\pcal,u\pmodels\neg\boxn n B$ and claim 1, we have 
$\neg\boxn n B\in u$. Then, since $\neg\boxn n B\to\boxn{n+1}\neg\boxn n B\in u$, 
we get $\boxn{n+1}\neg\boxn n B\in u$ and therefore, 
$\boxn {n+1}\big(\boxn{n+1}\neg\boxn n B \wedge \neg\boxn n B\big)\in u$. 
Hence by $u\R_{n+1} w$ we have $\boxn{n+1}\neg\boxn n B \wedge \neg\boxn n B\in w$.
This means $  \neg\boxn n B\in\axioms{\plgc w {n+1}}$, as desired.
\item Poly-necessitation: Let $\plgc w n\vdash B $. 
This means that $B\in\axioms{\plgc w n}$ and hence 
$B\wedge\boxn n B\in w$.
Then because $\boxn n B\to \boxn n \boxn n B\in w$, 
we get $\boxn n B\wedge\boxn n\boxn n B\in w$. Thus $\boxn n B\in\axioms{\plgc w n}$, as desired.
\item Poly-L\"ob's rule: Let $\plgc w n\vdash \boxn n B\to B$. 
Then by poly-necessitation, 
we get $\plgc w n\vdash \boxn n(\boxn n B\to B)$. On the other hand, 
obviously the L\"ob's axiom belongs to $\plgc w n$. 
Then, since $\plgc w n$ is closed under modus ponens, we get $\plgc w n\vdash \boxn n B$. 
Again modus ponens together with  $\plgc w n\vdash \boxn n B\to B$  
implies $\plgc w n \vdash B$, as desired.
\end{itemize}
%\\[2mm]
\textit{Claim 4.} There is a $w\in W$ such that 
$\pcal,w\npmodels A$  
while $\pcal,\sigma\pmodels B$ for every $B\in \Gamma$. 
\\[2mm]
\textit{Proof of the claim 4.}
First notice that together with previous claims, it finishes the proof of this Theorem.
So to prove this claim, we reason as follows. 
Since $\Gamma\nvdash_\GLP A$, there is a maximal consistent set 
$w\supseteq \Gamma\cup\GLP$ such that 
$A\nin w$. Then claim 1 implies $\pcal, w\npmodels A$ and $\pcal,w\pmodels  \Gamma$, 
as desired.
\end{proof}

\begin{question}\label{q-GLP}
The current poly-provability $\GLP$-models, are not having decidable validity. 
Also we have not constructed the provability models in a predicative way, 
as we did in the case of $\GL$ (\Cref{bases-gl}) or $\ILM$ (\Cref{bases-ilm}).
Our question here is to find a class of poly-provability models that are decidable,  
and for which, $\GLP$ is sound and complete. 
\end{question}
Probably, the main cause for non-predicativity of the poly-provability $\GLP$-models
defined in the completeness theorem \ref{GLP-completeness}, 
is that the accessibility relations are not converse well-founded. 
So probably, the crucial step would be to try to replace the logical property of 
poly-L\"ob's rule by its corresponding frame property, namely converse well-foundedness.

\section{Future works}

Given that provability models are well-aligned with logics with provability nature, 
it is very natural to ask for such semantics for the \textit{Logic of Proofs}  (${\sf LP}$)
(see \cite{Artemov}) and the justification logic of $\GL$, named ${\sf JGL}$ (see \cite{Shamkanov,GhariT}). The main idea of the justification logic is that, it explicitely
provides a proof-term (justification) for $\Box A$. In other words, it replaces $\Box A$
by $t:A$, meaning that $t$ is a proof/justification for $A$. 
\begin{itemize}
\item Define classes of provability-style models for which 
${\sf LP}$ and  ${\sf JGL}$ are sound and complete. 
\end{itemize}

The provability models, first used in \cite{PLHA} as a semantics for an intuitionistic 
modal logic $\iglh$ in the language $\lcalb$. This $\iglh$, is an extension of the 
intuitionistic $\GL$, annotated as $\igl$, with additional axioms 
of the form $\Box A\to\Box B$, for all  $\HA$-verifiable admissible rules 
\Ax{$A$}\UI{$B$}\DP of $\HA$. For example the axiom $\Box\neg\neg\Box A\to\Box\Box A$
belongs to $\iglh$, which corresponds to the Markov Rule:
\begin{center} 
\Ax{$\neg\neg\varphi$}\Ax{$\varphi\in\Sigma_1$}\BI{$\varphi$}\DP. 
\end{center}
On the other hands, we have intuitionistic modal logics in the language 
with two modal operators $\Box$ and $\lozenge$. Unlike the 
classical modal logic, $\lozenge$ is not definable via $\Box$ as $\neg\Box\neg$. 
Henceforth, we have dual axioms for $\lozenge$. See \cite{Simpson} for axiomatizations 
and Kripke models for such logics. Some easy questions like decidability, 
becomes challenging in such intuitionistic modal logics. For instance, 
it is due to a very recent result \cite{IS4} that we know Intuitionistic ${\sf S4}$ 
is decidable. In case of intuitionistic $\GL$, namely ${\sf IGL}$ 
(as introduced in \cite{Das}), even a Hibert-style 
axiomatization is missing! 
Then it is very natural question to ask about provability models of 
intuitionistic modal logics with two modal operators:
\begin{itemize} 
\item Define classes of provability models for which, various intuitionistic modal logics
in the language with both $\Box$ and $\lozenge$, are sound and complete. 
\end{itemize}

\section*{Acknowledgements}
The authors of this paper would like to express their 
deep gratitude to the following people for many fruitful
conversations that have contributed to this work:
Amir Akbar Tabatabai, Mohammad Ardeshir, David Fernández-Duque, Rosalie Iemhoff, 
Joost Joosten, Ian Shillito, Fedor Pakhomov, Albert Visser, and Ren-June Wang. 
In particular, they are especially 
indebted to Joost Joosten for kindly sending us the 
Master’s thesis
\cite{Mayaux} and explaining some ideas included in it.  
We also owe special thanks to Ian Shillito, 
whose question regarding the
strong completeness of $\GL$ for provability models motivated the inclusion of a subsection
on the strong completeness result.    
We also are grateful for the very pleasant and enlightening conversation with Amir Akbar
Tabatabai, in which he pointed out some interesting objections which cause improvement in 
the content of this paper.
Finally, we wish to  thank Fedor Pakhomov  for a thoughtful conversation, 
whose sharp critiques were invaluable in improving this research.

\section*{Funding}
The first author of this paper, acknowledges that this 
work is partially funded by FWO grant G0F8421N
and BOF grant BOF.STG.2022.0042.01.

\mojref

\begin{appendix}
\section[Propositional-Provability interpretations for Modal Logics]{Propositional-Provability interpretations for Modal Logics}
\label{ap-a}
In this appendix, we consider the following interpretation for the modal operator  $\Box$:
$\Box A$ is interpreted as $\sft\vdash A$,
for some propositional theory $\sft$ in the language $\lcalb$. Then we show that 
$\GL$ is sound for  $\sft$-provability iff 
$\sft$ includes $\GL$ and is closed under modus ponens 
and necessitation.

We first need to define precisely what we mean by such 
$\sft$-provability interpretation for an arbitrary 
$A\in\lcalb$.

A \textit{phrase} is a  formula of the shape $\bigwedge X\to\bigvee Y$, such that 
$X$ and $Y$ are finite sets of atomics or boxed formulas and $X\cap Y=\emptyset$. 
Notice that $\bot$ and $\top$ can not be members of $X\cup Y$.
Also for the uniqueness reasons for our later application, we assume that $\bigwedge X$ is the conjunction of members of $X$ with some fixed order on the set of atomics and boxed formulas. We do the same for the calculation of $\bigvee Y$. This avoids non-uniquenesses for our later applications.
Furthermore,  
we assume that $\bigwedge\emptyset:=\top$ and $\bigvee\emptyset:=\bot$.
For every given formula $A$, in classical logic, one may calculate a unique finite set 
$Z$ of phrases such that $\vdash A\lr\bigwedge Z$. 
Such set $Z$, is called 
the \textit{conjunctive normal form} of $A$. It is not difficult to observe that this conjunctive normal 
form is unique. 
This means that if $Z$ and $Z'$ are both sets of phrases
such that $\vdash\bigwedge Z\lr\bigwedge Z'$, 
then for every phrase $(\bigwedge X\to \bigvee Y)\in Z$
we also have $(\bigwedge X\to \bigvee Y)\in Z'$ and vice versa.

The $\sft$-provability interpretation, namely $A^\sft$, is defined as a statement in 
the metalanguage. We first define $A^\sft$ for a phrase $A=\bigwedge X\to\bigvee Y$ 
as follows. $A^\sft$ is defined to be true, if  the following holds:
\begin{itemize}
%\item There is some atomic $p\in X\cap Y$.
\item If for every $\Box E\in X$ we have $\sft\vdash E$, 
		then there is some $\Box F\in Y$ such that $\sft\vdash F$.
\end{itemize}
 
Then for the calculation of $A^\sft$ in general, we first calculate the 
conjunctive normal form of $A$, i.e.~ a set $Z$ 
of phrases, and then define 
$$A^\sft:=\bigwedge\{B^\sft: B\in Z\}.$$
Notice that the formula in the right hand side is an statement in the metalanguage.  
Finally we define 
the provability logic  $\PL(\sft)$ as the set of all formulas $A\in\lcalb$ such that 
$A^\sft$ is true.

\begin{lemma}\label{lem-cnf}
Two   classically equivalent formulas have equivalent 
$\sft$-provability interpretations. In other words, $\vdash A\lr B$  implies that  
$A^\sft$ holds, if and only if $B^\sft$ holds.
\end{lemma}
\begin{proof}
Since $A$ and $B$ are classically equivalent, they share same conjunctive normal form, and 
hence, they share equivalent 
$\sft$-provability interpretations.  
\end{proof}

\begin{theorem}[\textbf{Soundness.}]\label{GL-sound-Lint}
$\GL\subseteq\PL(\sft) $ iff  $\sft$ is a logic containing 
$\GL$ and closed under modus ponens and   necessitation.
\end{theorem}
\begin{proof}
First notice that if $\GL\subseteq\PL(\sft) $, then we have:
\begin{itemize}
\item $\sft$ contains all theorems of $\GL$: 
because for any theorem $A$ of $\GL$, we have 
	$\GL\vdash \Box A$ and hence we should have 
	$(\Box A)^\sft$, i.e.~$\sft\vdash A$.
\item $\sft$ is closed under modus ponens: because for $A:=\Box(B\to C)\to (\Box B\to\Box C)$ we have 
$\GL\vdash A$ and hence $A^\sft$ holds. 
\item $\sft$ is closed under necessitation: 
because for $A:=\Box B\to \Box \Box B$ we have $\GL\vdash A$ and hence $A^\sft$ holds. 
\end{itemize}
So it remains only to show the other way around. 
We use the axiomatization of $\GL$ without necessitation (see \cref{sec-gl}). 
It is straightforward to argue by induction on the proof $\GL\vdash A$ that $A^\sft$ holds.
Let us first treat the only inference rule, namely modus ponens. 
Assume that $\GL\vdash A$ and $\GL\vdash A\to B$, from which by modus ponens we inferred
$\GL\vdash B$. By induction hypothesis, we have the statements $A^\sft$ and $(A\to B)^\sft$.
Given that $(.)^\sft$ commutes with conjunction, we have also 
$(A\wedge (A\to B))^\sft$. Given that $A\wedge (A\to B)$ is classically equivalent 
to $A\wedge B$, by \Cref{lem-cnf} also $(A\wedge B)^\sft$ holds. Again, since 
$(.)^\sft$ commutes with conjunctions,  $B^\sft$ also holds. 
It remains  only to treat the axioms of $\GL$:
\begin{itemize}
\item $A$ is Classical tautology. Then the conjunctive normal form of $A$ 
is conjunction over empty set and hence $A^\sft$ is also the conjunction over an empty set,
thus holds in the metalanguage.
\item $A=\Box(B\to C)\to (\Box B\to\Box C)$. We have $A^\sft$ 
since $\sft$ is assumed to 	be closed under modus ponens.
\item $A=\Box B\to\Box \Box B$. 	We have $A^\sft$ since $\sft$
is assumed to be closed under necessitation.
\item $A=\Box(\Box B\to B)\to\Box B$. First notice that by L\"ob's axiom and necessitation 
and modus ponens, $\sft$ is also closed under L\"ob's rule. 
Then  $A^\sft$ holds because $\sft$ is closed under L\"ob's rule.
\item $A=\Box B$ and $B$ is any of above axioms. Then $A^\sft$ holds since 
 $\sft$ includes $\GL$.\qedhere
\end{itemize}
\end{proof}
As one should expect, we shouldn't have the completeness 
of $\GL$ for $\sft$-provability model for any theory $\sft$:
\begin{theorem}\label{GL-incompleteness}
For any theory $\sft$, we have $\GL\neq\PL(\sft)$. 
\end{theorem}
\begin{proof}
%Assume that  $\GL\subseteq \PL(\sft)$. 
We have either of following cases:
\begin{itemize}
\item $\sft\vdash \bot$: Then we have $\Box\bot\in \PL(\sft) \setminus\GL$.
\item $\sft\nvdash \bot$: Then we have $\neg\Box\bot\in\PL(\sft) \setminus \GL$.\qedhere
\end{itemize}
\end{proof}

\end{appendix}
\end{document}